\documentclass[12pt]{amsart}
\usepackage{amssymb}
\newtheorem{theorem}{Theorem}[section]
\newtheorem{lemma}[theorem]{Lemma}
\newtheorem{proposition}[theorem]{Proposition}
\newtheorem{corollary}[theorem]{Corollary}

\theoremstyle{definition}

\newtheorem{examp}[theorem]{Example}

\theoremstyle{remark}
\newtheorem{remark}[theorem]{Remark}

\numberwithin{equation}{section}

\newenvironment{example}{\begin{examp}\rm}{\diams\end{examp}}
\newcommand{\diams}{\unskip\nobreak\hfil\penalty50%
\hskip1em\hbox{}\nobreak\hfil%
$\diamondsuit$\parfillskip=0pt\finalhyphendemerits=0}

\newcommand{\bfind}[1]{\index{#1}{\bf #1}}
\newcommand{\n}{\par\noindent}

\newcommand{\sn}{\par\smallskip\noindent}
\newcommand{\mn}{\par\medskip\noindent}
\newcommand{\bn}{\par\bigskip\noindent}
\newcommand{\pars}{\par\smallskip}
\newcommand{\parm}{\par\medskip}
\newcommand{\parb}{\par\bigskip}
\newcommand{\isom}{\simeq}
\newcommand{\dec}{^d}
\newcommand{\Aut}{\mbox{\rm Aut}\,}

\newcommand{\ovl}[1]{\overline{#1}}
\newcommand{\ec}{\prec_{\exists}}
\newcommand{\sep}{^{\rm sep}}
\newcommand{\chara}{\mbox{\rm char}\,}
\newcommand{\trdeg}{\mbox{\rm trdeg}\,}
\newcommand{\Quot}{\mbox{\rm Quot}\,}
\newcommand{\Gal}{\mbox{\rm Gal}\,}
\newcommand{\rr}{\mbox{\rm rr}\,}
\newcommand{\bbox}[1]{\makebox(0,0){\rule[-2ex]{0ex}{5.4ex}#1}}

\newcommand{\cal}{\mathcal}
\newcommand{\N}{\mathbb N}
\newcommand{\Q}{\mathbb Q}
\newcommand{\R}{\mathbb R}
\newcommand{\Z}{\mathbb Z}

\newcommand{\F}{\mathbb F}
\newcommand{\Fp}{\F_p}
\begin{document}
\title{The defect}
\author{Franz-Viktor Kuhlmann}
\address{Department of Mathematics and
Statistics, University of Saskatchewan, 106 Wiggins Road,
Saskatoon, Saskatchewan, Canada S7N 5E6}
\email{fvk@math.usask.ca}
\thanks{This work was partially supported by a Canadian NSERC grant and
by a sabbatical grant from the University of Saskatchewan. I
thank Bernard Teissier, Olivier Piltant, Michael Temkin and Hagen
Knaf for helpful discussions and support. I thank the two referees for
their careful reading and useful comments.}
\date{13.\ 4.\ 2010}

\maketitle

\begin{abstract}
We give an introduction to the valuation theoretical
phenomenon of ``defect'', also known as ``ramification deficiency''. We
describe the role it plays in deep open
problems in positive characteristic: local uniformization (the local
form of resolution of singularities), the model theory of valued
fields, the structure theory of valued function fields. We give several
examples of algebraic extensions with non-trivial defect. We indicate
why Artin-Schreier defect extensions play a central role and describe a
way to classify them. Further, we give an overview of various results
about the defect that help to tame or avoid it, in particular
``stability'' theorems and theorems on ``henselian rationality'', and
show how they are applied. Finally, we include a list of open problems.
\end{abstract}

{\footnotesize\rm
\tableofcontents
}
%
%
%
%
\begin{section}{Valued fields}
Historically, there are three main origins of valued fields:
\sn
1) Number theory: Kurt Hensel introduced the fields $\Q_p$ of $p$-adic
numbers and proved the famous Hensel's Lemma (see below) for them. They
are defined as the completions of $\Q$ with respect to the (ultra)metric
induced by the $p$-adic valuations of $\Q$, similarly as the field of
reals, $\R$, is the completion of $\Q$ with respect to the usual metric
induced by the ordering on $\Q$.
\sn
2) Ordered fields: $\R$ is the maximal archimedean ordered field;
any ordered field properly containing $\R$ will have infinite
elements, that is, elements larger than all reals. Their inverses are
infinitesimals. The classes of magnitude, called archimedean classes,
give rise to a natural valuation. These valuations are important in
the theory of ordered fields and in real algebraic geometry.

In connection with ordered fields and their classes of magnitude, Hans
Hahn (\cite{[H]}) introduced an important class of valued fields, the
\textbf{(generalized) power series fields}. Take any field $K$ and any
ordered abelian group $G$. Let $K((G))$ (also denoted by $K((t^G))\,$)
be the set of all maps $\mu$ from $G$ to $K$ with well-ordered
\textbf{support} $\{g\in G\mid \mu(g)\ne 0\}$. One can visualize the
elements of $K((G))$ as formal power series $\sum_{g\in G}^{} c_g t^g$
for which the support $\{g\in G\mid c_g\ne 0\}$ is well-ordered. Using
this condition one shows that $K((G))$ is a field. Also, one uses it
to introduce the valuation:
\begin{equation}
v_t\,\sum_{g\in G}^{} c_g t^g\>=\>\min\{g\in G\mid c_g\ne 0\}
\end{equation}
(the minimum exists because the support is well-ordered). This valuation
is called the \textbf{canonical valuation} or \textbf{$t$-adic valuation}
of $K((G))$, and sometimes called the \textbf{minimum support valuation}.
Note that $v_t t=1$. For $G=\Z$, one obtains the field of formal Laurent
series $K((t))$.

\sn
3) Function fields: if $K$ is any field and $X$ an indeterminate, then
the rational function field $K(X)$ has a $p(X)$-adic valuation for every
irreducible polynomial $p(X)\in K[X]$, plus the $1/X$-adic valuation.
These valuations are trivial on $K$. As a valuation can be extended to
every extension field, these valuations together with the $p$-adic
valuations mentioned in 1) yield that a field admits a non-trivial
valuation as soon as it is not algebraic over a finite field. In
particular, all algebraic function fields over $K$ (i.e., finitely
generated field extensions of $K$ of transcendence degree $\geq 1$)
admit non-trivial valuations that are trivial on $K$. Such valued
function fields play a role in several areas of algebra and number
theory, some of which we will mention in this paper. Throughout,
\textbf{function field} will always mean {\it algebraic function field}.
%

\parm
If $K$ is a field with a valuation $v$, then we will denote its
value group by $vK$ and its residue field by $Kv$. For $a\in K$, its
value is $va$, and its residue is $av$. An extension of valued fields
is written as $(L'|L,v)$, meaning that $L'|L$ is a field extension,
$v$ is a valuation on $L'$ and $L$ is equipped with the restriction of
this valuation. Then there is a natural embedding of the value group
$vL$ in the value group $vL'$, and a natural embedding
of the residue field $Lv$ in the residue field $L'v$. If both embeddings
are onto (which we just express by writing $vL=vL'$ and $Lv=L'v$), then
the extension $(L'|L,v)$ is called {\bf immediate}. For $a\in L'$ we set
$v(a-L):=\{v(a-c)\mid c\in L\}$. The easy proof of the following
lemma is left to the reader:

\begin{lemma}                               \label{imexl}
The extension $(L'|L,v)$ is immediate if and only if for all $a\in L'$
there is $c\in L$ such that $v(a-c)>va$. If the extension $(L'|L,v)$
is immediate, then $v(a-L)$ has no maximal element and is an initial
segment of $vL$, that is, if $\alpha\in v(a-L)$ and $\alpha>\beta\in
vL$, then $\beta\in v(a-L)$.
\end{lemma}

If for each $a\in L'$ and every $\alpha\in vL'$ there is $c\in L$ such
that $v(a-c)>\alpha$, then we say that $(L,v)$ is \textbf{dense in}
$(L',v)$. If this holds, then the extension $(L'|L,v)$ is immediate. The
maximal extension in which $(L,v)$ is dense is its \textbf{completion}
$(L,v)^c$, which is unique up to isomorphism.

\pars
Every finite extension $L'$ of a valued field $(L,v)$
satisfies the \bfind{fundamental inequality} (cf.\ (17.5)
of \cite{[End]} or Theorem 19 on p.~55 of \cite{[Z--S]}):

\begin{equation}                             \label{fundineq}
n\>\geq\>\sum_{i=1}^{\rm g} {\rm e}_i {\rm f}_i
\end{equation}
where $n=[L':L]$ is the degree
of the extension, $v_1,\ldots,v_{\rm g}$ are the distinct extensions
of $v$ from $L$ to $L'$, ${\rm e}_i=(v_iL':vL)$ are the respective
\textbf{ramification indices} and ${\rm f}_i=[L'v_i:Lv]$ are the
respective \textbf{inertia degrees}. If ${\rm g}=1$ for every finite
extension $L'|L$ then $(L,v)$ is called \bfind{henselian}. This holds if
and only if $(L,v)$ satisfies \textbf{Hensel's Lemma}, that is, if $f$ is
a polynomial with coefficients in the valuation ring ${\cal O}$ of
$(L,v)$ and there is $b\in {\cal O}$ such that $vf(b)>0$ and $vf'(b)=0$,
then there is $a\in {\cal O}$ such that $f(a)=0$ and $v(b-a)>0$.

Every valued field $(L,v)$ admits a \bfind{henselization}, that is, a
minimal algebraic extension which is henselian (see
Section~\ref{sectram} below). All henselizations are isomorphic over
$L$, so we will frequently talk of {\it the} henselization of $(L,v)$,
denoted by $(L,v)^h$. The henselization becomes unique in absolute terms
once we fix an extension of the valuation $v$ from $L$ to its algebraic
closure. All henselizations are immediate separable-algebraic
extensions. If $(L',v)$ is a henselian extension field of $(L,v)$, then
a henselization of $(L,v)$ can be found inside of $(L',v)$.

\pars
For the basic facts of valuation theory, we refer the reader to
\cite[Appendix]{[Ax2]}, \cite{[End]}, \cite{[Eng--Pr]}, \cite{[Ri]},
\cite{[W]} and \cite{[Z--S]}. For ramification theory, we recommend
\cite{[End]}, \cite{[Eng--Pr]} and \cite{[Ne]}. For basic facts of model
theory, see \cite{[CK]}.

For a field $K$, $\tilde{K}$ will denote its algebraic closure and
$K\sep$ will denote its separable-algebraic closure. If $\chara K=p$,
then $K^{1/p^{\infty}}$ will denote its perfect hull. If we have
two subfields $K,L$ of a field $M$ (in our cases, we will usually have
the situation that $L\subset\tilde{K}$) then $K.L$ will denote the
smallest subfield of $M$ which contains both $K$ and $L$; it is called
the \textbf{field compositum of $K$ and $L$}.

\end{section}

%
%
\begin{section}{Two problems}                      \label{sect2p}
Let us look at two important problems that will lead us to considering
the phenomenon of defect:

%
%
\begin{subsection}{Elimination of ramification}
Given a valued function field
$(F|K,v)$, we want to find nice generators of $F$ over $K$. For
instance, if $F|K$ is separable then it is separably generated, that is,
there is a transcendence basis $T$ such that $F|K(T)$ is a finite
separable extension, hence simple. So we can write $F=K(T,a)$ with $a$
separable-algebraic over $K(T)$.

In the presence of the valuation $v$, we may want to ask for more. The
problem of \bfind{smooth local uniformization} is to find generators
$x_1,\ldots,x_n$ of $F|K$ in the valuation ring ${\cal O}$ of $v$ on $F$
such that the point $x_1v,\ldots,x_nv$ is smooth, that is, the Implicit
Function Theorem holds in this point. We say that $(F|K,v)$ is
\textbf{inertially generated} if there is a transcendence basis $T$ such
that $F$ lies in the \textbf{absolute inertia field} $K(T)^i$ (see
Section~\ref{sectram} for its definition). A connection between both
notions is given by Theorem~1.6 of \cite{[Kn-Ku2]}:

\begin{theorem}                             \label{lu=ig}
If $(F|K,v)$ admits smooth local uniformization, then it is
inertially generated.
\end{theorem}
If $(F|K,v)$ is inertially generated by the transcendence basis $T$,
then $vF=vK(T)$, and $Fv|K(T)v$ is separable. If this were not true, we
would say that $(F|K(T),v)$ is \textbf{ramified}. Let us consider an
example.

\begin{example}                             \label{ex1}
Suppose that $v$ is a discrete valuation on $F$ which is trivial on $K$
and such that $Fv|K$ is algebraic.
So there is an element $t\in F$ such that $vF=\Z vt=vK(t)$. Take the
henselization $F^h$ of $F$ with respect to some fixed extension of $v$
to the algebraic closure of $F$.

Assume that $\trdeg F|K =1$; then $F|K(t)$ is finite. Take $K(t)^h$ to
be the henselization of $K(t)$ within $F^h$. Then $F^h|K(t)^h$ is again
finite since $F^h=F.K(t)^h$ (cf.\ Theorem~\ref{hensfcomp} below). If
$\trdeg F|K >1$, we can take $T$ to be a transcendence basis of $F|K$
which contains $t$. Then again, $F^h|K(T)^h$ is finite, and $vF=vK(T)$.
But does that prove that $F|K$ is inertially generated? Well, if for
instance $K$ is algebraically closed, then it follows that
$Fv=K=K(T)v$, so $Fv|K(T)v$ is separable.
But ``inertially generated'' asks for more. In order to show that
$F|K$ is inertially generated in this particular case, we would have to
find $T$ such that $F\subseteq K(T)^h$, that is, the extension
$F^h|K(T)^h$ is trivial (see Section~\ref{sectram}).

Since $K(T)^h$ is henselian, there is only one extension of the given
valuation from $K(T)^h$ to $F^h$. By our choice of $T$, we have e$\>=
(vF^h:vK(T)^h)=(vF:vK(T))=1$, and if $K$ is algebraically closed, also
f$\>=(F^hv:K(T)^h v)=(Fv:K(T))v=1$. Hence equality holds in the
fundamental inequality (\ref{fundineq}) if and only if $F^h|K(T)^h$ is
trivial.
\end{example}

This example shows that it is important to know when the fundamental
inequality (\ref{fundineq}) is in fact an equality, or more precisely,
what the quotient $n/$ef is. A first and important answer is given by
the \bfind{Lemma of Ostrowski}. Assume that $(L'|L,v)$
is a finite extension and the extension of $v$ from $L$ to $L'$ is
unique. Then the Lemma of Ostrowski says that
\begin{equation}                            \label{LoO}
[L':L]\;=\; p^\nu \cdot (vL':vL)\cdot [L'v:Lv] \;\;\;\mbox{ with }
\nu\geq 0
\end{equation}
where $p$ is the \textbf{characteristic exponent} of $Lv$, that is,
$p=\chara Lv$ if this is positive, and $p=1$ otherwise. The Lemma of
Ostrowski can be proved using Tschirnhausen transformations (cf.\
\cite[Theoreme 2, p.~236]{[Ri]}). But it can also be deduced from
ramification theory, as we will point out in Section~\ref{sectram}
(see also \cite[Corollary to Theorem~25, p.~78]{[Z--S]}).

The factor ${\rm d}={\rm d}(L'|L,v)=p^\nu$ is called the \bfind{defect}
(or \textbf{ramification deficiency} as in \cite[p.~58]{[Z--S]}) of the
extension $(L'|L,v)$. If d$\>=1$, then we call $(L'|L,v)$ a
\bfind{defectless extension}; otherwise, we call it a \textbf{defect
extension}. Note that $(L'|L,v)$ is always defectless if $\chara Kv=0$.

We call $(L,v)$ a \bfind{defectless field}, \textbf{separably defectless
field} or \textbf{inseparably defectless field} if equality holds in the
fundamental inequality (\ref{fundineq}) for every finite, finite
separable or finite purely inseparable, respectively, extension
$L'$ of $L$. One can trace this back to the case of unique extensions of
the valuation; for the proof of the following theorem, see \cite{[Ku2]}
(a partial proof was already given in Theorem~18.2 of \cite{[End]}):
\begin{theorem}
A valued field $(L,v)$ is a defectless field if and only if its
henselization is. The same holds for ``separably defectless'' and
``inseparably defectless''.
\end{theorem}

Therefore, the Lemma of Ostrowski shows:
\begin{corollary}                               \label{lostr0}
Every valued field $(L,v)$ with $\chara Lv=0$ is a defectless field.
\end{corollary}

The defect is multiplicative in the following sense. Let $(L|K,v)$ and
$(M|L,v)$ be finite extensions. Assume that the extension of $v$ from
$K$ to $M$ is unique. Then the defect satisfies the following product
formula
\begin{equation}         \label{pf}
\mbox{\rm d}(M|K,v) = \mbox{\rm d}(M|L,v)\cdot\mbox{\rm d}(L|K,v)
\end{equation}
which is a consequence of the multiplicativity of the degree of field
extensions and of ramification index and inertia degree.
This formula implies:
\begin{lemma}                                         \label{md}
$(M|K,v)$ is defectless if and only if $(M|L,v)$ and $(L|K,v)$ are
defectless.
\end{lemma}
\begin{corollary}      \label{ed}
If $(L,v)$ is a defectless field and $(L',v)$ is a finite extension of
$(L,v)$, then $(L',v)$ is also a defectless field. Conversely, if there
exists a finite extension $(L',v)$ of $(L,v)$ such that $(L',v)$ is a
defectless field, the extension of $v$ from $L$ to $L'$ is unique,
and the extension $(L'|L,v)$ is defectless, then $(L,v)$ is a defectless
field. The same holds for ``separably defectless'' in the place of
``defectless'' if $L'|L$ is separable, and for ``inseparably
defectless'' if $L'|L$ is purely inseparable.
\end{corollary}

\pars
The situation of our Example~\ref{ex1} becomes more complicated when the
valuations are not discrete:

\begin{example}                             \label{ex2}
There are valued function fields $(F|K,v)$ of transcendence degree 2
with $v$ trivial on $K$ such that $vF$ is not finitely generated.
Already on a rational function field $K(x,y)$, the value group of a
valuation trivial on $K$ can be any subgroup of the rationals $\Q$ (see
Theorem 1.1 of \cite{[Ku6]} and the references given in that paper). In
such cases, if $F$ is not a rational function field, it is not easy to
find a transcendence basis $T$ such that $vF=vK(T)$. But even if we find
such a $T$, what do we know then about the extension $(F^h|K(T)^h,v)$?
For example, is it defectless?
\end{example}

An extension $(L'|L,v)$ of henselian fields is called
\textbf{unramified} if $vL'=vL$, $L'v|Lv$ is separable and every finite
subextension of $(L'|L,v)$ is defectless. Hence if $\chara Lv=0$, then
$(L'|L,v)$ is unramified already if $vL'=vL$. Note that our definition
of ``unramified'' is stronger than the definition in \cite[\$22]{[End]}
which does not require ``defectless''.

For a valued function field $(F|K,v)$, \bfind{elimination of
ramification} means to find a transcendence basis $T$ such that
$(F^h|K(T)^h,v)$ is unramified. According to Theorem~\ref{mue} in
Section~\ref{sectram} below, this is equivalent to $F$ lying in the
absolute inertia field $K(T)^i$. Hence, $(F|K,v)$ admits elimination of
ramification if and only if it is inertially generated.

If $\chara K=0$ and $v$ is trivial on $K$, then $(F|K,v)$ is always
inertially generated; this follows from Zariski's local uniformization
(\cite{[Z]}) by Theorem~\ref{lu=ig}. Since then $\chara Fv=\chara
Kv=\chara K=0$, Zariski did not have to deal with inseparable residue
field extensions and with defect. But if $\chara K>0$, then the
existence of defect makes the problem of local uniformization much
harder. This becomes visible in the approach to local uniformization
that is used in the papers \cite{[Kn-Ku1]} and \cite{[Kn-Ku2]}. Local
uniformization can be proved for Abhyankar places in positive
characteristic because the defect does not appear (\cite{[Kn-Ku1]}); we
will discuss this in more detail below. For other places
(\cite{[Kn-Ku2]}), the defect has to be ``killed'' by a finite extension
of the function field (``alteration'').

\end{subsection}


%
%
\begin{subsection}{Classification of valued fields up to elementary
equivalence}                                \label{subsee}
Value group and residue field are invariants of a valued field, that is,
two isomorphic valued fields have isomorphic value groups and isomorphic
residue fields. But two valued fields with the same value groups and
residue fields need not at all be isomorphic. For example, the valued
field $(\Fp(t),v_t)$ and $(\Fp((t)),v_t)$ both have value group $\Z$ and
residue field $\Fp$, but they are not isomorphic since $\Fp(t)$ is
countable and $\Fp((t))$ is not.

In situations where classification up to isomorphism fails,
classification up to elementary equivalence may still be possible. Two
algebraic structures are \textbf{elementarily equivalent} if they
satisfy the same elementary (first order) sentences. For example,
Abraham Robinson proved that all algebraically closed valued fields of
fixed characteristic are elementarily equivalent (cf.\
\cite[Theorem~4.3.12]{[Ro]}). James Ax and Simon Kochen and,
independently, Yuri Ershov proved that two henselian valued fields are
elementarily equivalent if their value groups are elementarily
equivalent and their residue fields are elementarily equivalent and of
characteristic 0 (cf.\ \cite{[AK1]} and \cite[Theorem~5.4.12]{[CK]})).
They also proved that all $p$-adically closed fields are elementarily
equivalent (cf.\ \cite[Theorem~2]{[AK2]}). Likewise, Alfred Tarski
proved that all real closed fields are elementarily equivalent (cf.\
\cite[Theorem~4.3.3]{[Ro]} or \cite[Theorem~5.4.4]{[CK]}). This remains
true if we consider non-archimedean real closed fields together with
their natural valuations (\cite{[CD]}). These facts (and the
corresponding model completeness results) all have important
applications in algebra (for instance, Nullstellens\"atze, Hilbert's
17th Problem, cf.\ \cite[Chap.~A 4, \$2]{[HB]}). So we would like to
know when classification up to elementary equivalence is possible for
more general classes of valued fields.

Two elementarily equivalent valued fields have elementarily equivalent
value groups and elementarily equivalent residue fields. When does the
converse hold? We mentioned already that the henselization is an
immediate extension. So the elementary properties of value group
and residue field do not determine whether a field is henselian or not.
But being henselian is an elementary property, expressed by a scheme of
elementary sentences, one each for all polynomials of degree $n$, where
$n$ runs through all natural numbers. In our above example,
$(\Fp((t)),v_t)$ is henselian, but $(\Fp(t),v_t)$ is not, so they are
not elementarily equivalent. We see that in order to have classification
up to elementary equivalence relative to value groups and residue
fields, our fields need to be (at least) henselian. But if the
characteristic of the residue fields is positive, then we may have
proper immediate algebraic extensions of henselian valued fields, as we
will see in the next section. So our fields need to be (at least)
\textbf{algebraically maximal},\index{algebraically maximal field} that
is, not admitting any proper immediate algebraic extensions.

Our fields even have to be defectless. Indeed, every valued field
$(K,v)$ admits a maximal immediate extension $(M,v)$. Then $(M,v)$ is
maximal and therefore henselian and defectless. Since $vK=vM$ and
$Kv=Mv$, we want that $(K,v)\equiv (M,v)$. The property ``henselian and
defectless field'' is elementary (cf.\ \cite[1.33]{[D]}, \cite{[Ku2]} or
the background information in \cite{[Ku9]}), so $(K,v)$ should be a
henselian defectless field.

\parm
If ${\cal L}$ is an elementary language and ${\cal A}\subset {\cal B}$
are ${\cal L}$-structures, then we will say that ${\cal A}$ is
\textbf{existentially closed in} ${\cal B}$ and write ${\cal A}\ec
{\cal B}$ if every existential sentence with parameters from ${\cal A}$
that holds in ${\cal B}$ also holds in ${\cal A}$. When we talk of
fields, then we use the language of rings ($\{+,-,\cdot,0,1\}$) or
fields (adding the unary function symbol ``$.^{-1}$''). When we talk of
valued fields, we augment this language by a unary relation symbol for
the valuation ring or a binary relation symbol for valuation
divisibility (``$vx\leq vy$''). For ordered abelian groups, we use the
language of groups augmented by a binary predicate (``$x<y$'') for the
ordering. For the meaning of ``existentially closed in'' in the settings
of fields, valued fields and ordered abelian groups, see
\cite[p.~183]{[Ku--Pr]}.

By model theoretical tools such as Robinson's Test, the
classification problem can be  transformed into the problem of finding
conditions which ensure that the following \bfind{Ax--Kochen--Ershov
Principle} holds:
\begin{equation}                            \label{AKE}
(K,v)\subseteq (L,v)\>\wedge\> vK\ec vL\>\wedge\> Kv\ec Lv
\;\;\;\Longrightarrow\;\;\; (K,v)\ec (L,v)\;.
\end{equation}
In order to prove that $(K,v)\ec (L,v)$, we first note that existential
sentences in $L$ only talk about finitely many elements of $L$, and
these generate a function field over $K$. So it suffices to show $(K,v)
\ec (F,v)$ for every function field $F$ over $K$ contained in $L$. One
tool to show that $(K,v)\ec (F,v)$ is to prove an embedding lemma: we
wish to construct an embedding of $(F,v)$ over $K$ in some ``big''
(highly saturated elementary) extension $(K^*,v^*)$ of $(K,v)$.
Existential sentences are preserved by embeddings and will then hold in
$(K^*,v^*)$ from where they can be pulled down to $(K,v)$. In order to
construct the embedding, we need to understand the algebraic structure
of $(F,v)$.

\begin{example}                             \label{exclhr}
Assume that $(K,v)$ is henselian (the same will then be true for
$(K^*,v^*)\,$) and that $(F|K,v)$ is an immediate extension of
transcendence degree 1. Pick an element $x\in F$ transcendental over
$K$. Even if we know how to embed $(K(x),v)$ in $(K^*,v^*)$, how can we
extend this embedding to $(F,v)$? Practically the only tool we have for
such extensions is Hensel's Lemma. So if $F\subset K(x)^h$, we can use
the universal property of henselizations (Theorem~\ref{hensuniqemb}
below) to extend the embedding to $K(x)^h$ and thus to $F$. If $F$ is
not a subfield of $K(x)^h$, we do not know what to do.
\end{example}

More generally, we have to deal with extensions which are not immediate,
but for which the conditions ``$vK\ec vL$'' and ``$Kv\ec Lv$'' hold. By
the saturation of $(K^*,v^*)$, they actually provide us with an
embedding of $vF$ over $vK$ in $v^*K^*$ and an embedding of $Fv$ over
$Kv$ in $K^*v^*$. Using Hensel's Lemma, they can be lifted to an
embedding of $(F,v)$ in $(K^*,v^*)$ if $(F,v)$ is inertially generated
with a transcendence basis $T$ such that $(K(T),v)$ can be embedded, as
we will discuss at the end of Section~\ref{subsst}. If $(F,v)$ is not
inertially generated, we are lost again. So we see that both of our
problems share the important approach of elimination of ramification.

\pars
Before we discuss the stated problems further, let us give several
examples of defect extensions, in order to meet the enemy we are dealing
with.

\end{subsection}
\end{section}

%
%
\begin{section}{Examples for non-trivial defect}    \label{sectexadef}
In this section, we shall give examples for extensions with defect $>1$.
There is one basic example which is quick at hand. It is due to
F.~K.~Schmidt.
\begin{example}                             \label{exampFKS}
We consider $\Fp((t))$ with its canonical valuation $v=v_t\,$.
Since $\F_p((t))|\F_p(t)$ has infinite transcendence degree, we can
choose some element $s\in\Fp((t))$ which
is transcendental over $\Fp(t)$. Since $(\Fp((t))|\Fp(t),v)$ is an
immediate extension, the same holds for $(\Fp(t,s)|\Fp(t),v)$ and thus
also for $(\Fp(t,s)|\Fp(t,s^p),v)$. The latter extension is purely
inseparable of degree $p$ (since $s,t$ are algebraically independent
over $\Fp\,$, the extension $\Fp(s)|\Fp(s^p)$ is linearly disjoint from
$\Fp(t,s^p)|\Fp(s^p)\,$). Hence, Theorem~\ref{allext} shows that
there is only one extension of the valuation $v$ from $\Fp(t,s^p)$ to
$\Fp(t,s)$. So we have ${\rm e}={\rm f}={\rm g}=1$ for this extension
and consequently, its defect is $p$.
\end{example}
\begin{remark}
This example is the easiest one used in commutative algebra to show that
the integral closure of a noetherian ring of dimension $1$ in a finite
extension of its quotient field need not be finitely generated.
\end{remark}

In some sense, the field $\Fp(t,s^p)$ is the smallest possible admitting
a defect extension. Indeed, a function field of transcendence degree $1$
over its prime field $\Fp$ is defectless under every valuation. More
generally, a valued function field of transcendence degree $1$ over a
subfield on which the valuation is trivial is always a defectless field;
this follows from Theorem~\ref{ai} below.

With respect to defects, discrete valuations are not too bad. The
following is easy to prove (cf.~\cite{[Ku2]}):

\begin{theorem}                             \label{vgZsdl}
Let $(K,v)$ be a discretely valued field, that is, with value group
$vK\isom\Z$. Then every finite separable extension is defectless. If in
addition $\chara K=0$, then $(K,v)$ is a defectless field.
\end{theorem}

A defect can appear ``out of nothing'' when a finite extension
is lifted through another finite extension:
\begin{example}                             \label{exampindef}
In the foregoing example, we can choose $s$ such that $vs>1=vt$. Now we
consider the extensions $(\Fp(t,s^p)|\Fp(t^p,s^p),v)$ and $(\Fp(t+s,s^p)
|\Fp(t^p,s^p),v)$ of degree $p$. Both are defectless: since
$v\Fp(t^p,s^p) =p\Z$ and $v(t+s)=vt=1$, the index of $v\Fp(t^p,s^p)$ in
$v\Fp(t,s^p)$ and in $v\Fp(t+s,s^p)$ must be (at least) $p$. But
$\Fp(t,s^p).\Fp(t+s, s^p) = \Fp(t,s)$, which shows that the defectless
extension $(\Fp(t,s^p)| \Fp(t^p,s^p),v)$ does not remain defectless if
lifted up to $\Fp(t+s,s^p)$ (and vice versa).
\end{example}

We can derive from Example~\ref{exampFKS} an example of a defect
extension of henselian fields.
\begin{example}                             \label{exampFKSh}
We consider again the immediate extension $(\Fp(t,s)| \Fp(t,s^p),v)$ of
Example~\ref{exampFKS}. We take the henselization $(\Fp(t,s),v)^h$ of
$(\Fp(t,s),v)$ in $\Fp((t))$ and the henselization $(\Fp(t,s^p),v)^h$ of
$(\Fp(t,s^p),v)$ in $(\Fp(t,s),v)^h$. We find that $(\Fp(t,s),v)^h|
(\Fp(t,s^p),v)^h$ is again a purely inseparable extension of degree $p$.
Indeed, the purely inseparable extension $\Fp(t,s)| \Fp(t, s^p)$ is
linearly disjoint from the separable extension
$\Fp(t,s^p)^h|\Fp(t,s^p)$, and by virtue of Theorem~\ref{hensfcomp},
$\Fp(t,s)^h= \Fp(t,s).\Fp(t,s^p)^h$. Also for this extension we have
that ${\rm e}={\rm f}={\rm g}=1$ and again, the defect is $p$. Note that
by Theorem~\ref{vgZsdl}, a proper immediate extension over a field like
$(\Fp(t,s^p),v)^h$ can only be purely inseparable.
\end{example}

\parm
The next example is easily found by considering the purely inseparable
extension $\tilde{K}|K\sep$. In comparison to the last example, the
involved fields are ``much bigger'', for instance, they do not have
value group $\Z$ anymore.

\begin{example}                             \label{exaperfhull}
Let $K$ be a field which is not perfect. Then the extension $\tilde{K}|
K\sep$ is non-trivial. For every non-trivial valuation $v$ on
$\tilde{K}$, the value groups $v\tilde{K}$ and $vK\sep$ are both equal
to the divisible hull $\widetilde{vK}$ of $vK$, and the residue fields
$\tilde{K}v$ and $K\sep v$ are both equal to the algebraic closure of
$Kv$ (cf.\ Lemma~2.16 of \cite{[Ku6]}). Consequently, $(\tilde{K}|
K\sep,v)$ is an immediate extension. Since the extension of $v$ from
$K\sep$ to $\tilde{K}$ is unique (cf.\ Theorem~\ref{allext} below), we
find that the defect of every finite subextension is equal to its
degree.
\end{example}

Note that the separable-algebraically closed field $K\sep$ is henselian
for every valuation. Hence, our example shows:

\begin{theorem}                             \label{hinsimm}
There are henselian valued fields of positive characteristic which admit
proper purely inseparable immediate extensions. Hence, the property
``henselian'' does not imply the property ``algebraically maximal''.
\end{theorem}

We can refine the previous example as follows. Let $p>0$ be the
characteristic of the residue field $Kv$. An \bfind{Artin-Schreier
extension} of $K$ is an extension of degree $p$ generated over $K$ by a
root of a polynomial $X^p-X-c$ with $c\in K$. An extension of degree $p$
of a field of characteristic $p$ is a Galois extension if and only if it
is an Artin-Schreier extension. A field $K$ is \textbf{Artin-Schreier
closed} if it does not admit Artin-Schreier extensions.

\begin{example}
In order that every purely inseparable extension of the valued field
$(K,v)$ be immediate, it suffices that $vK$ be $p$-divisible and $Kv$ be
perfect. But these conditions are already satisfied for every
non-trivially valued Artin-Schreier closed field $K$ (see
Corollary~2.17 of \cite{[Ku6]}). Hence, the perfect hull of every
non-trivially valued Artin-Schreier closed field is an immediate
extension.
\end{example}

Until now, we have only presented purely inseparable defect extensions.
But our last example can give an idea of how to produce a separable
defect extension by interchanging the role of purely inseparable
extensions and Artin-Schreier extensions.

\begin{example}                             \label{exampiASe}
Let $(K,v)$ be a valued field of characteristic $p>0$ whose value group
is not $p$-divisible. Let $c\in K$ such that $vc<0$ is not divisible by
$p$. Let $a$ be a root of the Artin-Schreier polynomial $X^p-X-c$. Then
$va=vc/p$ and $[K(a):K]=p= (vK(a):vK)$. The fundamental inequality shows
that $K(a)v=Kv$ and that the extension of $v$ from $K$ to $K(a)$ is
unique. By Theorem~\ref{allext} below the further extension to
$K(a)^{1/p^{\infty}}= K^{1/p^{\infty}}(a)$ is unique. It follows that
the extension of $v$ from $K^{1/p^{\infty}}$ to $K^{1/p^{\infty}}(a)$ is
unique. On the other hand, $[K^{1/p^{\infty}} (a):K^{1/p^{\infty}}]=p$
since the separable extension $K(a)|K$ is linearly disjoint from
$K^{1/p^{\infty}}|K$. The value group $vK^{1/p^{\infty}}(a)$ is the
$p$-divisible hull of $vK(a)= vK+\Z va$. Since $pva\in vK$, this is the
same as the $p$-divisible hull of $vK$, which in turn is equal to
$vK^{1/p^{\infty}}$. The residue field of $K^{1/p^{\infty}}(a)$ is the
perfect hull of $K(a)v=Kv$. Hence it is equal to the residue field of
$K^{1/p^{\infty}}$. It follows that the extension $(K^{1/p^{\infty}}(a)|
K^{1/p^{\infty}},v)$ is immediate and that its defect is $p$, like its
degree.

Similarly, one can start with a valued field $(K,v)$ of characteristic
$p>0$ whose residue field is not perfect. In this case, the
Artin-Schreier extension $K(a)|K$ is constructed as in the proof of
Lemma~2.13 of \cite{[Ku6]}. We leave the details to the reader.
\end{example}

In the previous example, we can always choose $(K,v)$ to be henselian
(since passing to the henselization does not change value group and
residue field). Then all constructed extensions of $(K,v)$ are also
henselian, since they are algebraic extensions (cf.\
Theorem~\ref{hensfcomp} below). Hence, our example shows:

\begin{theorem}                             \label{hsepimm}
There are henselian valued fields of positive characteristic which admit
immediate Artin-Schreier defect extensions.
\end{theorem}

If the perfect hull of a given valued field $(K,v)$ is not an immediate
extension, then $vK$ is not $p$-divisible or $Kv$ is not perfect,
and we can apply the procedure of our above example. This shows:
\begin{theorem}
If the perfect hull of a given valued field of positive characteristic
is not an immediate extension, then it admits an immediate
Artin-Schreier extension.
\end{theorem}

An important special case of Example~\ref{exampiASe} is the following:
\begin{example}                             \label{fse}
We choose $(K,v)$ to be $(\Fp(t),v_t)$ or $(\Fp((t)),v_t)$ or any
intermediate field, and set $L:=K(t^{1/p^i}\mid i\in\N)$, the perfect
hull of $K$. By Theorem~\ref{allext} below, $v=v_t$ has a unique
extension to $L$. In all cases, $L$ can be viewed as a subfield of the
power series field $\Fp((\Q))$. The power series
\begin{equation}
\vartheta\>:=\>\sum_{i=1}^{\infty} t^{-1/p^i}\;\in\> \Fp((\Q))
\end{equation}
is a root of the Artin-Schreier polynomial
\[X^p-X-\frac{1}{t}\]
because
\begin{eqnarray*}
\vartheta^p-\vartheta-\frac{1}{t} & = & \sum_{i=1}^{\infty}
t^{-1/p^{i-1}}\,-\,\sum_{i=1}^{\infty} t^{-1/p^i}\,-\,t^{-1} \\
& = & \sum_{i=0}^{\infty} t^{-1/p^{i}}\,-\,\sum_{i=1}^{\infty}
t^{-1/p^i}\,-\,t^{-1} \>=\> 0\;.
\end{eqnarray*}
By Example~\ref{exampiASe}, the extension $L(\vartheta)|L$ is an immediate
Artin-Schreier defect extension. The above power series expansion for
$\vartheta$ was presented by Shreeram Abhyankar in \cite{[Ab1]}. It
became famous since it shows that there are elements algebraic over
$\Fp(t)$ with a power series expansion in which the exponents do not
have a common denominator. This in turn shows that Puiseux series fields
in positive characteristic are in general not algebraically closed (see
also \cite{[Ke],[Ku4]}). With $p=2$, the above was also used by
Irving Kaplansky in \cite[Section~5]{[Ka]} for the construction of an
example that shows that if his ``hypothesis A'' (see
\cite[Section~3]{[Ka]}) is violated, then the maximal immediate
extension of a valued field may not be unique up to isomorphism. See
also \cite{[Ku--Pa--Ro]} for more information on this subject.

\pars
Let us compute $v(\vartheta-L)$. For the partial sums
\begin{equation}
\vartheta_k\>:=\>\sum_{i=1}^k t^{-1/p^i}\;\in\> L
\end{equation}
we see that $v(\vartheta-\vartheta_k)=-1/p^{k+1}<0$. Assume that there
is $c\in L$ such that $v(\vartheta-c)>-1/p^k$ for all $k$. Then
$v(c-\vartheta_k)=\min\{v(\vartheta-c),v(\vartheta-\vartheta_k)\}
=-1/p^{k+1}$ for all $k$. On the other hand, there is some $k$ such that
$c\in K(t^{-1/p},\ldots,t^{-1/p^k})=K(t^{-1/p^k})$. But this
contradicts the fact that $v(c-t^{-1/p}-\ldots-t^{-1/p^k})=v(c-
\vartheta_k)= -1/p^{k+1}\notin vK(t^{-1/p^k})$. This proves that the
values $-1/p^k$ are cofinal in $v(\vartheta-L)$. Since $vL$ is a
subgroup of the rationals, this shows that the least upper bound of
$v(\vartheta-L)$ in $vL$ is the element 0. As $v(\vartheta-L)$ is an
initial segment of $vL$ by Lemma~\ref{imexl}, we conclude that
$v(\vartheta-L)=(vL)^{<0}$. It follows that $(L(\vartheta)|L,v)$ is
immediate without $(L,v)$ being dense in $(L(\vartheta),v)$.
\end{example}

A version of this example with $(K,v)=(\widetilde{\Fp}((t)),v_t)$ was
given by S.~K.~Khanduja in \cite{[Kh]} as a counterexample to
Proposition~$2'$ on p.~425 of \cite{[Ax1]}. That proposition states that
if $(K,v)$ is a perfect henselian valued field of rank 1 and $a\in
\tilde{K}\setminus K$, then there is $c\in K$ such that
\[
v(a-c)\>\geq\>\min\{v(a-a')\mid a'\ne a \mbox{ conjugate to $a$ over
$K$}\}\;.
\]
But for $a=\vartheta$ in the previous example, we have that $a-a'\in
\Fp$ so that the right hand side is $0$, whereas $v(\vartheta-c)<0$ for
all $c$ in the perfect hull $L$ of $\Fp((t))$. The same holds if we take
$L$ to be the perfect hull of $K=\widetilde{\Fp}((t))$. In fact, it is
Corollary~2 to Lemma~6 on p.~424 in \cite{[Ax1]} which is in error; it
is stated without proof in the paper.

\pars
In a slightly different form, the previous example was already given by
Alexander Ostrowski in \cite{[Os]}, Section~57:

\begin{example}
Ostrowski takes $(K,v)= (\Fp(t), v_t)$, but works with the
polynomial $X^p-tX-1$ in the place of the Artin-Schreier polynomial
$X^p-X-1/t$. After an extension of $K$ of degree $p-1$, it also can be
transformed into an Artin-Schreier polynomial. Indeed, if we take $b$ to
be an element which satisfies $b^{p-1}=t$, then replacing $X$ by $bX$
and dividing by $b^p$ will transform $X^p-tX-1$ into the polynomial
$X^p-X-1/b^p$. Now we replace $X$ by $X+1/b$. Since we are working in
characteristic $p$, this transforms $X^p-X-1/b^p$ into $X^p-X-1/b$.
(This sort of transformation plays a crucial role in the proofs of
Theorem~\ref{ai} and Theorem~\ref{hr} as well as in Abhyankar's and
Epp's work.). Now we see that the Artin-Schreier polynomial $X^p-X-1/b$
plays the same role as $X^p-X-1/t$. Indeed, $vb=\frac{1}{p-1}$ and it
follows that $(v\Fp(b):v\Fp(t))=p-1=[\Fp(b):\Fp(t)]$, so that $v\Fp(b)=
\Z\frac{1} {p-1}$. In this value group, $vb$ is not divisible by $p$.
\end{example}

\pars
Interchanging the role of purely inseparable and Artin-Schreier
extensions in Example~\ref{fse}, we obtain:

\begin{example}                             \label{fsei}
We proceed as in Example~\ref{fse}, but replace $t^{-1/p^i}$ by
$a_i$, where we define $a_1$ to be a root of the Artin-Schreier
polynomial $X^p-X-1/t$ and $a_{i+1}$ to be a root of the Artin-Schreier
polynomial $X^p-X+a_i\,$. Now we choose $\eta$ such that
$\eta^p=1/t$. Note that also in this case, $a_1,\ldots,a_i\in
K(a_i)$ for every $i$, because $a_i= a_{i+1}^p -a_{i+1}$ for every
$i$. By induction on $i$, we again deduce that $va_1= -1/p$ and
$va_i= -1/p^i$ for every $i$. We set $L:=K(a_i\mid i\in \N)$,
that is, $L|K$ is an infinite tower of Artin-Schreier extensions. By our
construction, $vL$ is $p$-divisible and $Lv=\Fp$ is perfect. On the
other hand, for every purely inseparable extension $L'|L$ the group
$vL'/vL$ is a $p$-group and the extension $L'v|Lv$ is purely
inseparable. This fact shows that $(L(\eta)|L,v)$ is an immediate
extension.

In order to compute $v(\eta-L)$, we set
\begin{equation}
\eta_k\>:=\>\sum_{i=1}^{k} a_i\;\in\> L\;.
\end{equation}
Bearing in mind that $a_{i+1}^p=a_{i+1}-a_i$ and $a_1^p=
a_1+1/t$ for $i\geq 1$, we compute
\begin{eqnarray*}
(\eta-\eta_k)^p & = & \eta^p-\eta_k^p \;=\;
\frac{1}{t}\>-\>\sum_{i=1}^{k} a_i^p\;=\;
\frac{1}{t}\>-\left(\>\sum_{i=1}^{k} a_i\>-\>
\sum_{i=1}^{k-1} a_i\>+\>\frac{1}{t}\right)\\
 & = & a_k\;.
\end{eqnarray*}
It follows that $v(\eta-\eta_k)=\frac{va_k}{p}= -1/p^{k+1}\,$.
The same argument as in Example~\ref{fse} now shows that again,
$v(\eta-L)=(vL)^{<0}$.
\end{example}

\pars
We can develop Examples~\ref{fse} and~\ref{fsei} a bit further in order
to treat complete fields.
\begin{example}                             \label{complnotdl}
Take one of the immediate extensions $(L(\vartheta)|L,v)$ of
Example~\ref{fse} and set $\zeta=\vartheta$, or take one of the
immediate extensions $(L(\eta)|L,v)$ of Example~\ref{fsei} and set
$\zeta=\eta$. Consider the completion $(L,v)^c=(L^c,v)$ of $(L,v)$.
Since every finite extension of a complete valued field is again
complete, $(L^c(\zeta),v) =(L(\zeta).L^c,v)$ is the completion of
$(L(\zeta),v)$ for every extension of the valuation $v$ from $(L^c,v)$
to $L(\zeta).L^c$. Consequently, the extension $(L^c(\zeta)|L(\zeta),v)$
and thus also the extension $(L^c(\zeta)|L,v)$ is immediate. It follows
that $(L^c(\zeta)| L^c,v)$ is immediate. On the other hand, this
extension is non-trivial since $v(\zeta -L)=(vL)^{<0}$ shows that
$\zeta\notin L^c$.
\end{example}

A valued field is \textbf{maximal}\index{maximal field} if it does not
admit any proper immediate extension. All power series fields are
maximal. A valued field is said to have \bfind{rank 1} if its value
group is archimedean, i.e., a subgroup of the reals. Every complete
discretely valued field of rank 1 is maximal. Every complete valued
field of rank 1 is henselian (but this is not true in general in higher
ranks). The previous example proves:

\begin{theorem}
There are complete fields of rank 1 which admit immediate
separable-algebraic and immediate purely inseparable extensions.
Consequently, not every complete field of rank 1 is maximal.
\end{theorem}

In Example~\ref{fsei} we constructed an immediate purely inseparable
extension not contained in the completion of the field. Such extensions
can be transformed into immediate Artin-Schreier defect
extensions:

\begin{example}                             \label{transf}
In the situation of Example~\ref{fsei}, extend $v$ from $L(\eta)$ to
$\tilde{L}$. Take $d\in L$ with $vd\geq 1/p$, and $\vartheta_0$ a root
of the polynomial $X^p-dX-1/t$. It follows that
\[-1 \>=\> v\frac{1}{t} \>=\> v(\vartheta_0^p-d\vartheta_0)\>\geq\>
\min\{v\vartheta_0^p, vd\vartheta_0\} \>=\>
\min\{pv\vartheta_0, vd+v\vartheta_0\}\]
which shows that we must have $v\vartheta_0<0$. But then
\[pv\vartheta_0 \><\>v\vartheta_0\><\> vd+v\vartheta_0\;,\]
so \[v\vartheta_0 \>=\> -\frac{1}{p}\;.\]
We compute:
\[pv(\vartheta_0-\eta) \>=\>
v(\vartheta_0-\eta)^p \>=\> v(\vartheta_0^p-\eta^p) \>=\>
v(d\vartheta_0+1/t-1/t) \>=\> vd + v\vartheta_0 \>\geq\> 0\;.\]
Hence $v(\vartheta_0-\eta)\geq 0$, and thus for all $c\in L$,
\[v(\vartheta_0-c) \>=\> \min\{v(\vartheta_0-\eta),v(\eta-c)\}
\>=\> v(\eta-c)\;.\]
In particular, $v(\vartheta_0-L)=v(\eta-L)=(vL)^{<0}$. The extension
$(L(\vartheta_0)|L,v)$ is immediate and has defect $p\,$; however, this
is not quite as easy to show as it has been before. To make things
easier, we choose $(K,v)$ to be henselian, so that also $(L,v)$, being
an algebraic extension, is henselian. So there is only one extension of
$v$ from $L$ to $L(\vartheta_0)$. Since $v(\vartheta_0-c)<0$ for all $c
\in L$, we have that $\vartheta_0\notin L$. We also choose $d=b^{p-1}$
for some $b\in L$. Then we will see below that $L(\vartheta_0)|L$ is an
Artin-Schreier extension. If it were not immediate, then e$\>=p$ or
f$\>=p$. In the first case, we can choose some $a\in L(\vartheta_0)$
such that $0,va,\ldots, (p-1)va$ are representatives of the distinct
cosets of $vL(\vartheta_0)$ modulo $vL$. Then $1,a,\ldots,a^{p-1}$ are
$L$-linearly independent and thus form an $L$-basis of $L(\vartheta_0)$.
Writing $\vartheta_0=c_0+ c_1a+\ldots+c_{p-1}a^{p-1}$, we find that
$v(\eta-c_0)=v(\vartheta_0-c_0)=\min\{vc_1+va,\ldots,vc_{p-1}+(p-1)va\}
\notin vL$ as the values $vc_1+va,\ldots,vc_{p-1}+(p-1)va$ lie in
distinct cosets modulo $vL$. But this is a contradiction. In the second
case, f$\>=p$, one chooses $a\in L(\vartheta_0)$ such that $1,av,\ldots,
(av)^{p-1}$ form a basis of $L(\vartheta_0)v|Lv$, and derives a
contradiction in a similar way. (Using this method one actually proves
that an extension $(L(\zeta)|L,v)$ of degree $p$ with unique extension
of the valuation is immediate if and only if $v(\zeta-L)$ has no maximal
element.)

Now consider the polynomial $X^p-dX-1/t=X^p-b^{p-1}X-1/t$ and set
$X=bY$. Then $X^p-dX-1/t=b^pY^p-b^pY-1/t$, and dividing by $b^p$ we
obtain the polynomial $Y^p-Y-1/b^pt$ which admits $\vartheta_0/b$ as a
root. So we see that $(L(\vartheta_0)|L,v)$ is in fact an immediate
Artin-Schreier defect extension. But in comparison with
Example~\ref{fsei}, something is different:
\begin{eqnarray*}
v(\frac{\vartheta_0}{b}\,-\,L) & = & \{v(\frac{\vartheta_0}{b}\,-\,c)
\mid c\in L\} \>=\>\{v(\frac{\vartheta_0}{b}\,-\,\frac{c}{b})\mid
c\in L\} \\
& = & \{v(\vartheta_0-c)-vb\mid c\in L\} \>=\>\{\alpha\in vL\mid
\alpha<vb\}\;,
\end{eqnarray*}
where $vb>0$.
\end{example}

\pars
A similar idea can be used to turn the defect extension of
Example~\ref{exampFKS} into a separable extension. However, in the
previous example we made use of the fact that $\eta$ was not an element
of the completion of $(L,v)$, that is, $v(\eta-L)$ was bounded
from above. We use a ``dirty trick'' to first transform the extension of
Example~\ref{exampFKS} to an extension whose generator does not lie in
the completion of the base field.

\begin{example}                             \label{transfdt}
Taking the extension $(\Fp(t,s)|\Fp(t,s^p),v)$ as in
Example~\ref{exampFKS}, we adjoin a new transcendental element $z$ to
$\Fp(t,s)$ and extend the valuation $v$ in such a way that $vs\gg vt$,
that is, $v\Fp(t,s,z)$ is the lexicographic product $\Z\times\Z$. The
extension $(\Fp(t,s,z)|\Fp(t,s^p,z),v)$ is still purely inseparable and
immediate, but now $s$ does not lie anymore in the completion
$\Fp(t,s^p)((z))$ of $\Fp(t,s^p,z)$. In fact, $v(s-\Fp(t,s^p,z))=
\{\alpha\in v\Fp(t,s^p,z)\mid \exists n\in\N: nvt\geq\alpha\}$ is
bounded from above by $vz$.

Taking $\vartheta_0$ to be a root of the polynomial $X^p-z^{p-1}X-s^p$
we obtain that $v(\vartheta_0-c)=v(s-c)$ for all $c\in \Fp(t,s^p,z)$ and
that the Artin-Schreier extension $(\Fp(t,\vartheta_0,z)| \Fp(t,s^p,z)
,v)$ is immediate with defect $p$. We leave the proof as an exercise to
the reader. Note that one can pass to the henselizations of all fields
involved, cf.\ Example~\ref{exampFKSh}.
\end{example}

\parb
The interplay of Artin-Schreier extensions and radical extensions that
we have used in the last examples can also be transferred to the mixed
characteristic case. There are infinite algebraic extensions of $\Q_p$
which admit immediate Artin-Schreier defect extensions. To present an
example, we need a lemma which shows that there is some
\textbf{quasi-additivity} in the mixed characteristic case.

\begin{lemma}                               \label{quasadd}
Let $(K,v)$ be a valued field of characteristic 0 and
residue characteristic $p>0$, and with valuation ring ${\cal O}$.
Further, let $c_1,\ldots,c_n$ be elements in $K$ of value $\geq
-\frac{vp}{p}$. Then
\[(c_1+\ldots+c_n)^p\>\equiv\> c_1^p+\ldots+c_n^p
\pmod{\cal O}\;.\]
\end{lemma}
\begin{proof}
Every product of $p$ many $c_i$'s has value $\geq -vp$. In view of the
fact that every binomial coefficient
{\small\rm $\left(\begin{array}{c}  p\\ i \end{array}\right)$}
is divisible by $p$ for $1\leq i\leq p-1$, we find that $(c_1+c_2)^p
\equiv c_1^p+c_2^p \pmod{\cal O}$. Now the assertion follows by
induction on $n$.
\end{proof}

\begin{example}                             \label{Qpdefext}
We choose $(K,v)$ to be $(\Q,v_p)$ or $(\Q_p,v_p)$ or any intermediate
field. Note that we write $vp=1$. We construct an algebraic extension
$(L,v)$ of $(K,v)$ with a $p$-divisible value group as follows. By
induction, we choose elements $a_i$ in the algebraic closure of $K$
such that $a_1^p=1/p$ and $a_{i+1}^p=a_i\,$. Then $va_1=-1/p$ and
$va_i=-1/p^i$ for every $i$. Hence, the field $L:=K(a_i\mid i\in\N)$
must have $p$-divisible value group under any extension of $v$ from $K$
to $L$. Note that $a_1,\ldots, a_i\in K(a_i)$ for every $i$. Since
$(vK(a_{i+1}): vK(a_i)) =p$, the fundamental inequality shows that
$K(a_{i+1})v=K(a_i)v$ and that the extension of $v$ is unique, for every
$i$. Hence, $Lv=\Q_p v=\Fp$ and the extension of $v$ from $K$ to $L$ is
unique.

Now we let $\vartheta$ be a root of
$X^p-X-1/p$. It follows that $v\vartheta=-1/p$.
%
%
We define $b_i:=\vartheta-a_1-\ldots-a_i\,$. By construction, $va_i\geq
-1/p$ for all $i$. It follows that also $vb_i\geq -1/p$ for all $i$.
With the help of the foregoing lemma, and bearing in mind that
$a_{i+1}^p=a_i$ and $a_1^p=1/p$, we compute
\begin{eqnarray*}
0 & = & \vartheta^p-\vartheta-\frac{1}{p} \;=\;
(b_i+a_1+\ldots+a_i)^p-(b_i+a_1+\ldots+a_i)-1/p\\
 & \equiv & b_i^p-b_i +a_1^p+\ldots+a_i^p-a_1-\ldots-a_i-1/p
\;=\; b_i^p-b_i-a_i \pmod{\cal O}\;.
\end{eqnarray*}
Since $va_i<0$, we have that $vb_i=\frac{1}{p}va_i=-1/p^{i+1}\,$. Hence,
$(vK(\vartheta,a_i):vK(a_i))=p= [K(\vartheta, a_i):K(a_i)]$ and
$K(\vartheta,a_i)v=K(a_i)v=\Fp$ for every $i$. If $[L(\vartheta):L]<p$,
then there would exist some $i$ such that $[K(\vartheta,a_i):K(a_i)]<p$.
But we have just shown that this is not the case. Similarly, if
$vL(\vartheta)$ would contain an element that does not lie in the
$p$-divisible hull of $\Z=vK$, or if $L(\vartheta)v$ would be a proper
extension of $\Fp\,$, then the same would already hold for
$K(\vartheta,a_i)$ for some $i$. But we have shown that this is not the
case. Hence, $(L(\vartheta)|L,v)$ is an Artin-Schreier defect extension.

For the partial sums $\vartheta_k=\sum_{i=1}^{k} a_i$ we obtain
$v(\vartheta-\vartheta_k)= vb_k= -1/p^{k+1}$, and the same argument as
in Example~\ref{fse} shows again that $v(\vartheta-L)=(vL)^{<0}$.
\end{example}

From this example we can derive a special case which was given by
Ostrowski in \cite[Section 39]{[Os]} (see also \cite[Chap.~VI \$8,
exercise~2]{[B]}.
\begin{example}
In the last example, we take $K=\Q_2$. Then $(L(\sqrt{3})|L,v)$ is an
immediate extension of degree $2$. Indeed, this is nothing else than the
Artin-Schreier extension that we have constructed. If one substitutes
$Y=1-2X$ in the minimal polynomial $Y^2-3$ of $\sqrt{3}$ and then
divides by $4$, one obtains the Artin-Schreier polynomial $X^p-X-1/2$.

This is Ostrowski's original example. A slightly different version was
presented by Paulo Ribenboim (cf.\ Exemple 2 of Chapter G, p.~246): The
extension $(L(\sqrt{-1})|L,v)$ is immediate. Indeed, the minimal
polynomial $Y^2+1$ corresponds to the Artin-Schreier polynomial
$X^p-X+1/2$ which does the same job as $X^p-X-1/2$.
\end{example}

As in the equal characteristic case, we can interchange the role of
radical extensions and Artin-Schreier extensions:
\begin{example}
We proceed as in Example~\ref{Qpdefext}, with the only difference that
we define $a_1$ to be a root of the Artin-Schreier polynomial
$X^p-X-1/p$ and $a_{i+1}$ to be a root of the Artin-Schreier polynomial
$X^p-X+a_i\,$, and that we choose $\eta$ such that $\eta^p=
1/p$. Note that also in this case, $a_1,\ldots,a_i\in K(a_i)$ for
every $i$, because $a_i= a_{i+1}^p -a_{i+1}$ for every $i$. By induction
on $i$, we again deduce that $va_1= -1/p$ and that $va_i= -1/p^i$ for
every $i$. As before, we define $b_i:=a-a_1- \ldots-a_i\,$. Using
Lemma~\ref{quasadd} and bearing in mind that $a_{i+1}^p=a_{i+1}-a_i$ and
$a_1^p= a_1+1/p$, we compute
\begin{eqnarray*}
0 & = & \vartheta^p-\frac{1}{p} \;=\; (b_i+a_1+\ldots+a_i)^p-1/p\\
 & \equiv & b_i^p+a_1^p+\ldots+a_i^p-1/p
\;=\; b_i^p+a_i \pmod{\cal O}\;.
\end{eqnarray*}
It follows that $v(b_i^p+a_i)\geq 0>va_i\,$. Consequently, $vb_i^p=
-1/p^{i+1}\,$, that is, $vb_i=\frac{1}{p}va_i =va_{i+1}\,$. As before,
we set $L:=K(a_i\mid i\in \N)$. Now the same arguments as in
Example~\ref{Qpdefext} show that $(L(\vartheta)|L,v)$ is an immediate
extension with $v(\vartheta-L)=(vL)^{<0}$.
\end{example}

\pars
It can happen that it takes just a finite defect extension to make a
field defectless and even maximal. The following example is due to
Masuyoshi Nagata (\cite[Appendix, Example~(E3.1),
pp.~206-207]{[Na]}):
\begin{example}                             \label{exampnagata}
We take a field $k$ of
characteristic $p$ and such that $[k:k^p]$ is infinite, e.g., $k=
\Fp(t_i|i\in\N)$ where the $t_i$ are algebraically independent elements
over $\Fp$. Taking $t$ to be another transcendental element over $k$ we
consider the power series fields $k((t))$ and $k^p((t))= k^p((t^p))(t) =
k((t))^p (t)$. Since $[k:k^p]$ is not finite, we have that
$k((t))|k^p((t)).k$ is a non-trivial immediate purely inseparable
algebraic extension. In fact, a power series in $k((t))$ is an element
of $k^p((t)).k$ if and only if its coefficients generate a finite
extension of $k^p$. Since $k^p((t)).k$ contains $k((t))^p$, this
extension is generated by a set $X=\{x_i\mid i\in I\}\subset k((t))$
such that $x_i^p\in k^p((t)).k$ for every $i\in I$. Assuming this set to
be minimal, or in other words, the $x_i$ to be $p$-independent over
$k^p((t)).k$, we pick some element $x\in X$ and put $K:=k^p((t)).
k(X\setminus\{x\})$. Then $k((t))|K$ is a purely inseparable extension
of degree $p$. Moreover, it is an immediate extension; in fact, $k((t))$
is the completion of $K$. As an algebraic extension of $k^p((t))$, $K$
is henselian.
\end{example}

This example proves:
\begin{theorem}
There is a henselian discretely valued field $(K,v)$ of characteristic
$p$ admitting a finite immediate purely inseparable extension
$(L|K,v)$ of degree $p$ such that $(L,v)$ is complete, hence maximal
and thus defectless.
\end{theorem}

\parm
For the conclusion of this section, we shall give an example which is
due to Francoise Delon (cf.\ \cite{[D]}, Exemple~1.51). It shows
that an algebraically maximal field is not necessarily a defectless
field, and that a finite extension of an algebraically maximal field is
not necessarily again algebraically maximal.
%

\begin{example}                             \label{exdelon}
We consider $\Fp((t))$ with its $t$-adic valuation $v=v_t\,$. We choose
elements $x,y\in \Fp((t))$ which are algebraically independent over
$\Fp(t)$. We set $L:=\Fp(t,x,y)$ and define
\[s:=x^p+ty^p\;\mbox{\ \ and\ \ }\; K:=\Fp(t,s)\;.\]
Then $s$ is transcendental over $\F_p(t)$ and therefore, $K$ has
$p$-degree $2$, that is, $[K:K^p]=p^2$. We take $F$ to be the relative
algebraic closure of $K$ in $\Fp((t))$. Since the elements $1,t^{1/p},
\ldots, t^{(p-1)/p}$ are linearly independent over $\Fp((t))$, the same
holds over $F$. Hence, the elements $1,t,\ldots,t^{p-1}$ are linearly
independent over $F^p$. Now if $F$ had $p$-degree $1$, then $s$ could be
written in a unique way as an $F^p$-linear combination of $1,t,\ldots,
t^{p-1}$. But this is not possible since $s=x^p+ty^p$ and $x,y$ are
transcendental over $F$. Hence, the $p$-degree of $F$ is still $2$ (as
it cannot increase through algebraic extensions). On the other hand,
$vF=v\Fp((t))=\Z$ and $Fv=\Fp((t))v=\F_p\,$, hence $(vF:pvF)=p$ and
$[Fv:Fv^p]=1$. Now Theorem~\ref{charinsdl} shows that $(F,v)$ is not
inseparably defectless. Again from Theorem~\ref{charinsdl}, we infer
that $F^{1/p}=F(t^{1/p},s^{1/p})$ must be an extension of $F$ with
non-trivial defect. So $F$ is not a defectless field.

On the other hand, $\Fp((t))$ is the completion of $F$ since it is
already the completion of $\F_p(t)\subseteq F$. This shows that
$\Fp((t))$ is the unique maximal immediate extension of $F$ (up to
valuation preserving isomorphism over $F$). If $F$ would admit a proper
immediate algebraic extension $F'$, then a maximal immediate extension
of $F'$ would also be a maximal immediate extension of $F$ and would
thus be isomorphic over $F$ to $\Fp((t))$. But we have chosen $F$ to be
relatively algebraically closed in $\Fp((t))$. This proves that $(F,v)$
must be algebraically maximal.

As $(F,v)$ is algebraically maximal, the extension $F^{1/p}|F$ cannot
be immediate. Therefore, the defect of $F^{1/p}|F$ implies that both
$F^{1/p}|F(s^{1/p})$ and $F^{1/p}|F(t^{1/p})$ must be non-trivial
immediate extensions. Consequently, $F(s^{1/p})$ and
$F(t^{1/p})$ are not algebraically maximal.

Let us add to Delon's example by analyzing the situation in more detail
and proving that $F$ is the henselization of $K$ and thus a separable
extension of $K$. To this end, we first prove that $K$ is relatively
algebraically closed in $L$. Take $b\in L$ algebraic over $K$. The
element $b^p$ is algebraic over $K$ and lies in $L^p=\Fp(t^p,x^p,y^p)$
and thus also in $K(x)= \Fp(t,x,y^p)$. Since $x$ is transcendental over
$K$, $K$ is relatively algebraically closed in $K(x)$ and thus, $b^p\in
K$. Consequently, $b\in K^{1/p}= \Fp(t^{1/p},s^{1/p})$. Write
\[
b\;=\; r_0+r_1s^{\frac{1}{p}}+\ldots+r_{p-1}s^{\frac{p-1}{p}}
\mbox{ \ \ \ with \ } r_i\in \Fp(t^{1/p},s)=K(t^{1/p})\;.
\]
By the definition of $s$,
\[
b\;=\; r_0 \,+\, r_1x +\ldots+ r_{p-1}x^{p-1}
\,+\ldots+\, t^{1/p}r_1y +\ldots+ t^{(p-1)/p}r_{p-1}y^{p-1}
\]
(in the middle, we have omitted the summands in which both $x$ and $y$
appear). Since $x,y$ are algebraically independent over $\Fp$, the
$p$-degree of $\F_p(x,y)$ is $2$, and the elements $x^i y^j$, $0\leq
i<p$, $0\leq j<p$, form a basis of $\Fp(x,y)|\Fp(x^p,y^p)$. Since $t$
and $t^{1/p}$ are transcendental over $\Fp(x^p,y^p)$, we know that
$\Fp(x,y)|\Fp(x^p,y^p)$ is linearly disjoint from $\Fp(t,x^p,y^p)|
\Fp(x^p,y^p)$ and from $\Fp(t^{1/p},x^p,y^p)|\Fp(x^p,y^p)$. This shows
that the elements $x^i y^j$ also form a basis of $L|\Fp(t,x^p,y^p)$ and
are still $\Fp(t^{1/p},x^p,y^p)$-linearly independent. Hence, $b$ can
also be written as a linear combination of these elements with
coefficients in $\Fp(t,x^p,y^p)$, and this must coincide with the above
$\Fp(t^{1/p},x^p,y^p)$-linear combination which represents $b$. That is,
all coefficients $r_i\,$ {\it and\/} $\,t^{i/p}r_i$, $1\leq i<p$, are in
$\Fp(t,x^p,y^p)$. Since $t^{i/p}\notin\Fp(t,x^p,y^p)$, this is
impossible unless they are zero. It follows that $b=r_0\in K(t^{1/p})$.
Assume that $b\notin K$. Then $[K(b):K]=p$ and thus, $K(b)=K(t^{1/p})$
since also $[K(t^{1/p}):K]=p$. But then $t^{1/p}\in K(b)\subset L$, a
contradiction. This proves that $K$ is relatively algebraically closed
in $L$.

On the other hand, $t^{1/p}= y^{-1}(s^{1/p}-x)\in L(s^{1/p})$. Hence,
$L.K^{1/p}=L(t^{1/p},s^{1/p})=L(s^{1/p})$ and $[L.K^{1/p}:L]=
[L(s^{1/p}):L] \leq p<p^2=[K^{1/p}:K]$, that is, $L|K$ is not linearly
disjoint from $K^{1/p}|K$ and thus not separable. Although being
finitely generated, $L|K$ is consequently not separably generated;
in particular, it is not a rational function field. At the same time, we
have seen that $K(s^{1/p})$ admits a non-trivial purely inseparable
algebraic extension in $L(s^{1/p})$ (namely, $K^{1/p}$). In contrast,
$K(s^{1/p})$ and $L$ are $K$-linearly disjoint because $s^{1/p}
\notin L$.

Let us prove even more: if $K_1|K$ is any proper inseparable algebraic
extension, then $t^{1/p}\in L.K_1\,$. Take such an extension $K_1|K$.
Then there is some separable-algebraic subextension $K_2|K$ and an
element $a\in K_1 \setminus K_2$ such that $a^p\in K_2\,$. Since $K_2|K$
is separable and $K$ is relatively algebraically closed in $L$, we see
that $K_2$ is relatively algebraically closed in $L_2:=L.K_2\,$. Hence,
$a\notin L_2$ and therefore, $[L_2(a):L_2]=p$. On the other hand,
$K_2^{1/p}= K^{1/p} .K_2$ and thus, $L_2.K_2^{1/p}=L_2.K^{1/p}=
L.K^{1/p}.K_2\,$. Consequently, $[L.K^{1/p}:L]=p$ implies that
$[L_2.K_2^{1/p}:L_2]=[L.K^{1/p}.K_2:L.K_2]\leq p$. Since $a\in K_2^{1/p}
\subset L_2.K_2^{1/p}$ and $[L_2(a):L_2]=p$, it follows that
$L_2.K_2^{1/p}=L_2(a)$. We obtain:
\[t^{1/p}\in K^{1/p}\>\subseteq\>K_2^{1/p}\>\subseteq\>L_2.K_2^{1/p}
\>=\>L_2(a)\>\subseteq\>L.K_1\;.\]

If $F|K$ were inseparable, then $t^{1/p}\in L.F$, which contradicts
the fact that $L.F\subseteq \F_p((t))$. This proves that $F|K$ is
separable. Since $F$ is relatively closed in the henselian field
$\F_p((t))$, it is itself henselian and thus contains the henselization
$K^h$ of $K$. Now $\F_p((t))$ is the completion of $K^h$ since it is
already the completion of $\F_p(t)\subseteq K^h$. Since a henselian
field is relatively separable-algebraically closed in its completion
(cf.\ \cite{[W]}, Theorem~32.19), it follows that $F=K^h$.

Note that the maximal immediate extension $\Fp((t))$ of $F$ is not a
separable extension since its subextension $L.F|F$ is not linearly
disjoint from $K^{1/p}|K$.
\end{example}

This example proves:
\begin{theorem}
There are algebraically maximal fields which are not inseparably
defectless. Hence, ``algebraically maximal'' does not imply
``defectless''. There are algebraically maximal fields admitting a
finite purely inseparable extension which is not an algebraically
maximal field.
\end{theorem}

\end{section}

%
%
\begin{section}{Absolute ramification theory}       \label{sectram}
Assume that $L|K$ is an algebraic extension, not
necessarily finite, and that $v$ is a {\it non-trivial} valuation on $K$.
Throuthout this section, we fix an arbitrary extension of $v$ to the
algebraic closure $\tilde{K}$ of $K$, which we again denote by $v$. Then
for every $\sigma\in\Aut(\tilde{K}|K)$, the map
\begin{equation}
v\sigma=v\circ\sigma:\;L\ni a \>\mapsto\>v(\sigma a)\in v\tilde{K}
\end{equation}
is a valuation of $L$ which extends $v$. All extensions of $v$ from $K$
to $L$ are conjugate:
\begin{theorem}                             \label{allext}
The set of all extensions of $v$ from $K$ to $L$ is
\[\{v\sigma\mid \sigma \mbox{ an embedding of $L$ in
$\tilde{K}$ over $K$}\}\;.\]
In particular, a valuation on $K$ has a unique extension to every purely
inseparable field extension of $K$.
\end{theorem}

\pars
We will now give a quick introduction to \textbf{absolute ramification
theory}, that is, the ramification theory of the extension $\tilde{K}|K$
with respect to a given valuation $v$ on $\tilde{K}$ with valuation
ring ${\cal O}_{\tilde{K}}$. For a corresponding quick introduction to
general ramification theory, see \cite{[Ku4]}.

We define distinguished subgroups of the \textbf{absolute Galois group}
$G:=\Gal K:=\Aut(\tilde{K}|K)=\Aut(K\sep|K)$ of $K$ (with respect to the
fixed extension of $v$ to $\tilde{K}$). The
subgroup

\begin{equation}                            \label{decgrp1}
G\dec :=\{\sigma\in G\mid v\sigma=v\mbox{ on } \tilde{K}\}
\end{equation}
is called the \textbf{absolute decomposition group} {\bf of} $(K,v)$
(w.r.t.\ $(\tilde{K},v)$). Further, the \textbf{absolute inertia group}
(w.r.t.\ $(\tilde{K},v)$) is defined to be
\begin{equation}                            \label{ingrp}
G^i := \{\sigma\in G\mid \forall x\in {\cal O}_{\tilde{K}}:
v(\sigma x - x)>0\}\;,
\end{equation}
and the \textbf{absolute ramification group} (w.r.t.\ $(\tilde{K},v)$) is
\begin{equation}                            \label{ramgrp}
G^r:=G^r (L|K,v) :=\{\sigma\in G\mid
\forall x\in {\cal O}_{\tilde{K}}\setminus\{0\}: v(\sigma x - x)>vx\}\;.
\end{equation}
The fixed fields $K\dec$, $K^i$ and $K^r$ of $G\dec$, $G^i$ and $G^r$,
respectively, in $K\sep$ are called the \textbf{absolute decomposition
field}, \textbf{absolute inertia field} and \textbf{absolute ramification
field} {\bf of} $(K,v)$ (with respect to the given extension of $v$
to $\tilde{K}$).

\begin{remark}
In contrast to the classical definition used by other authors, we take
decomposition field, inertia field and ramification field to be the
fixed fields of the respective groups {\it in the separable-algebraic
closure of $K$}. The reason for this will become clear later.
\end{remark}

By our definition, $K\dec$, $K^i$ and $K^r$ are separable-algebraic
extensions of $K$, and $K\sep|K^r$, $K\sep|K^i$, $K\sep|K\dec$ are (not
necessarily finite) Galois extensions. Further,
\begin{equation}
1\subset G^r\subset G^i\subset G\dec\subset G\;\mbox{ and thus, }\;
K\sep\supset K^r\supset K^i\supset K\dec\supset K\;.
\end{equation}
(For the inclusion $G^i\subset G\dec$ note that $vx\geq 0$ and
$v(\sigma x-x)>0$ implies that $v\sigma x\geq 0$.)

\begin{theorem}
$G^i$ and $G^r$ are normal subgroups of $G\dec$, and $G^r$ is a normal
subgroup of $G^i$. Therefore, $K^i|K\dec$, $K^r|K\dec$ and
$K^r|K^i$ are (not necessarily finite) Galois extensions.
\end{theorem}

First, we consider the decomposition field $K\dec$. In some sense, it
represents all extensions of $v$ from $K$ to $\tilde{K}$.
\begin{theorem}                             \label{Z}
a) \ $v\sigma=v\tau$ on $\tilde{K}$ if and only if $\sigma\tau^{-1}$ is
trivial on $K\dec$.\n
b) \ $v\sigma=v$ on $K\dec$ if and only if $\sigma$ is trivial
on $K\dec$.\n
c) \ The extension of $v$ from $K\dec$ to $\tilde{K}$ is unique.\n
d) \ The extension $(K\dec|K,v)$ is immediate.
\end{theorem}
WARNING: It is in general not true that $v\sigma\ne v\tau$ holds already
on $K\dec$ if it holds on $\tilde{K}$.

Assertions a) and b) are easy consequences of the definition of $G\dec$.
Part c) follows from b) by Theorem~\ref{allext}. For d), there is a
simple proof using a trick mentioned by James Ax in
\cite[Appendix]{[Ax2]}; see also \cite[Theorem~22, p.~70 and
Theorem~23, p.~71]{[Ax2]} and \cite{[Eng--Pr]}.

\parm
Now we turn to the inertia field $K^i$. Let ${\cal M}_{\tilde{K}}$
denote the valuation ideal of $v$ on $\tilde{K}$ (the unique maximal
ideal of ${\cal O}_{\tilde{K}}$). For every $\sigma\in G\dec$ we have
that $\sigma {\cal O}_{\tilde{K}}={\cal O}_{\tilde{K}}$, and it follows
that $\sigma {\cal M}_{\tilde{K}} ={\cal M}_{\tilde{K}}$. Hence, every
such $\sigma$ induces an automorphism $\ovl{\sigma}$ of
${\cal O}_{\tilde{K}}/{\cal M}_{\tilde{K}}=\tilde{K}v=\widetilde{Kv}$
which satisfies $\ovl{\sigma}(av)= (\sigma a)v$. Since $\sigma$ fixes
$K$, it follows that $\ovl{\sigma}$ fixes $Kv$.

\begin{lemma}
The map
\begin{equation}                            \label{redofauto}
G\dec \ni\sigma\;\mapsto\;\ovl{\sigma}\in\Gal Kv
\end{equation}
is a group homomorphism.
\end{lemma}

\begin{theorem}                             \label{T}
a) \ The homomorphism (\ref{redofauto}) is onto and induces an
isomorphism
\begin{equation}
\Aut(K^i|K\dec)\>=\>G\dec/G^i\>\isom\> \Aut(K^i v|K\dec v)\;.
\end{equation}
b) \ For every finite subextension $F|K\dec$ of $K^i|K\dec$,
\begin{equation}
[F:K\dec]\>=\>[Fv:K\dec v]\;.
\end{equation}
c) \ We have that $vK^i=vK\dec =vK$. Further, $K^iv$ is the
separable closure of $Kv$, and therefore,
\begin{equation}
\Aut(K^i v|K\dec v)\>=\>\Gal Kv\;.
\end{equation}
\end{theorem}
If $F|K\dec$ is normal, then b) is an easy consequence of a). From
this, the general assertion of b) follows by passing from $F$ to the
normal hull of the extension $F|K\dec$ and then using the
multiplicativity of the extension degree. c) follows from b) by use of
the fundamental inequality.

\parm
We set $p:=\chara Kv$ if this is positive, and $p:=1$ if $\chara Kv=0$.
Given any abelian group $\Delta$, the \textbf{$p'$-divisible hull of
$\Delta$} is defined to be the subgroup $\{\alpha\in\tilde{\Delta}\mid
\exists n\in\N:\;(p,n) =1\, \wedge\,n\alpha\in\Delta\}$ of all elements
in the divisible hull $\tilde{\Delta}$ of $\Delta$ whose order modulo
$\Delta$ is prime to $p$.

\begin{theorem}                             \label{V}
a) \ There is an isomorphism
\begin{equation}
\Aut (K^r|K^i)\>=\> G^i/G^r\>\isom\>
\mbox{\rm Hom}\left(vK^r/vK^i\,,\,(K^rv)^\times\right)\;,
\end{equation}
where the character group on the right hand side is the full character
group of the abelian group $vK^r/vK^i$. Since this group is abelian,
$K^r|K^i$ is an abelian Galois extension.
\n
b) \ For every finite subextension $F|K^i$ of $K^r|K^i$,
\begin{equation}
[F:K^i]\>=\>(vF:vK^i)\;.
\end{equation}
c) \ $K^r v=K^i v$, and $vK^r$ is the $p'$-divisible hull of $vK$.
\end{theorem}
Part b) follows from part a) since for a finite extension $F|K^i$, the
group $vF/vK^i$ is finite and thus there exists an isomorphism of
$vF/vK^i$ onto its full character group. The equality $K^r v=
K^i v$ follows from b) by the fundamental inequality. The second
assertion of part c) follows from the next theorem and the fact that the
order of all elements in $(K^i v)^\times$ and thus also of all
elements in $\mbox{\rm Hom} \left(vK^r/vK^i\,,\,(K^i v)^\times\right)$
is prime to $p$.

\begin{theorem}                             \label{S}
The ramification group $G^r$ is a $p$-group, hence
$K\sep|K^r$ is a $p$-extension. Further, $v\tilde{K}/vK^r$ is a
$p$-group, and the residue field extension $\tilde{K}v|K^r v$ is
purely inseparable. If $\chara Kv=0$, then $K^r=K\sep=\tilde{K}$.
\end{theorem}
We note:
\begin{lemma}                              \label{RemASe}
Every $p$-extension is a tower of Galois extensions of degree $p$. In
characteristic $p$, all of them are Artin-Schreier extensions.
\end{lemma}

\pars
From Theorem~\ref{S} it follows that there is a canonical isomorphism
\begin{equation}                            \label{cg}
\mbox{\rm Hom}\left(vK^r/vK^i\,,\,(K^i v)^\times\right) \>\isom\>
\mbox{\rm Hom}\left(v\tilde{K}/vK\,,\,(\tilde{K}v)^\times\right)\;.
\end{equation}

\pars
The following theorem will be very useful for our purposes:
\begin{theorem}                             \label{liftZTV}
If $K'|K$ is algebraic, then the absolute decomposition field of
$(K',v)$ is $K\dec.K'$, its absolute inertia field is $K^i.K'$, and its
absolute ramification field is $K^r.K'\,$.
\end{theorem}

From part c) of Theorem~\ref{Z} we infer that the extension of $v$ from
$K\dec$ to $\tilde{K}$ is unique. On the other hand, if $L$ is any
extension field of $K$ within $K\dec$, then by Theorem~\ref{liftZTV},
$K\dec=L\dec$. Thus, if $L\ne K\dec$, then it follows from part b) of
Theorem~\ref{Z} that there are at least two distinct extensions of $v$
from $L$ to $K\dec$ and thus also to $\tilde{K}=\tilde{L}$. This proves
that the absolute decomposition field $K\dec$ is a minimal algebraic
extension of $K$ admitting a unique extension of $v$ to its algebraic
closure. So it is the minimal algebraic extension of $K$ which is
henselian. We call it the \textbf{henselization of $(K,v)$ in
$(\tilde{K},v)$}.\index{henselization} Instead of $K\dec$, we also write
$K^h$. A valued field is henselian if and only if it is equal to its
henselization. Henselizations have the following universal property:

\begin{theorem}                             \label{hensuniqemb}
Let $(K,v)$ be an arbitrary valued field and $(L,v)$ a henselian
extension field of $(K,v)$. Then there is a unique embedding of
$(K^h,v)$ in $(L,v)$ over $K$.
\end{theorem}

From part d) of Theorem~\ref{Z}, we obtain another very important
property of the henselization:
\begin{theorem}                             \label{immhens}
The henselization $(K^h,v)$ is an immediate extension of $(K,v)$.
\end{theorem}

\begin{corollary}                           \label{maxhens}
Every algebraically maximal and every maximal valued field is henselian.
In particular, $(K((t)),v_t)$ is henselian.
\end{corollary}

We employ Theorem~\ref{liftZTV} again to obtain:
\begin{theorem}                             \label{hensfcomp}
If $K'|K$ is an algebraic extension, then the henselization of
$K'$ is $K'.K^h\,$. Every algebraic extension of a henselian
field is again henselian.
\end{theorem}

\pars
In conjunction with Theorems~\ref{T} and~\ref{V}, Theorem~\ref{liftZTV}
is also used to prove that there are no defects between $K^h$ and $K^r$:

\begin{theorem}
Take a finite extension $K_2|K_1$ such that $K^h\subseteq K_1\subseteq
K_2\subseteq K^r$. Then $(K_2|K_1,v)$ is defectless.
\end{theorem}
\begin{proof}
Since $K_3:=K_2\cap K_1^i$ is a finite subextension of $K_1^i|K_1$, we
have by parts b) and c) of Theorem~\ref{T} that $[K_3:K_1]=[K_3v:K_1v]$
and $vK_3=vK_1$.
Since
$K_1^i|K_1$ is Galois, $K_2$ is linearly disjoint from $K_1^i$ over
$K_3$. That is, $[K_2.K_1^i:K_1^i]=[K_2:K_3]$. By Theorem~\ref{liftZTV},
$K_2\subseteq K^r = K_1.K^r=K_1^r$, so also $K_2.K_1^i|K_1^i$ is a
finite subextension of $K_1^r|K_1^i$. By part b) of Theorem~\ref{V}, we
thus have $[K_2.K_1^i:K_1^i]= (v(K_2.K_1^i):vK_1^i)$. By
Theorem~\ref{liftZTV}, $K_1^i=K_1^i.K_3=K_3^i$ and $K_2.K_1^i=K_2^i$, so
by part c) of Theorem~\ref{T}, $vK_1^i=vK_3^i=vK_3$ and
$v(K_2.K_1^i)=vK_2^i=vK_2$. Therefore,
\[[K_2:K_3] \>=\> [K_2.K_1^i:K_1^i] \>=\> (v(K_2.K_1^i):vK_1^i) \>=\>
(vK_2:vK_3) \>=\> (vK_2:vK_1)\]
Putting everything together, we obtain
\begin{eqnarray*}
[K_2:K_1] & = & [K_2:K_3][K_3:K_1] \>=\> (vK_2:vK_1)[K_3v:K_1v]\\
 & \leq & (vK_2:vK_1)[K_2v:K_1v] \>\leq\> [K_2:K_1]\;,
\end{eqnarray*}
so that equality must hold everywhere, which shows that
$(K_2|K_1,v)$ is defectless.
\end{proof}

An algebraic extension of $K^h$ is called \textbf{purely wild}\index{purely
wild extension} if it is linearly disjoint from $K^r$ over $K^h$. The
following theorem has been proved by Matthias Pank (see Theorem~4.3 and
Proposition~4.5 of \cite{[Ku--Pa--Ro]}):

\begin{theorem}                             \label{pank}
Every maximal purely wild extension $W$ of $K^h$ satisfies $W.K^r=
\tilde{K}$ and hence is a field complement of $K^r$ in $\tilde{K}$.
Moreover, $W^r=\tilde{W}$, $vW$ is the $p$-divisible hull of $vK$,
and $Wv$ is the perfect hull of $Kv$.
\end{theorem}

\begin{lemma}                             \label{dram}
If $(L|K^h,v)$ is a finite extension, then its
defect is equal to the defect of $(L.K^r|K^r,v)$.
\end{lemma}
\begin{proof}
We put $L_0 := L\cap K^r$. We have $L.K^r=L^r$ and $L_0^r=K^r$ by
Theorem~\ref{liftZTV}. Since $K^r|K^h$ is normal, $L$ is linearly
disjoint from $K^r=L_0^r$ over $L_0\,$, and $(L|L_0,v)$ is thus a purely
wild extension.

As a finite subextension of $(K^r|K^h,v)$, the extension $(L_0|K^h,v)$
is defectless. Hence by the multiplicativity of the defect (\ref{pf}),
\begin{equation}                            \label{d=d}
\mbox{\rm d}(L|K^h,v) = \mbox{\rm d}(L|L_0,v)\;.
\end{equation}

It remains to show $\mbox{\rm d}(L|L_0,v)=\mbox{\rm d}(L.K^r|K^r,v)$.
Since $L|L_0$ is linearly disjoint from $K^r|L_0\,$, we have
\begin{equation}                                  \label{gr0}
[L^r:L_0^r]\>=\> [L.K^r:L_0^r]\>=\> [L:L_0]\;.
\end{equation}
Since $L|L_0$ is purely wild, $vL/vL_0$ is a $p$-group and
$Lv|L_0v$ is purely inseparable. On the other hand, by Theorem~\ref{V},
\pars
$vL^r$ is the $p'$-divisible hull of $vL$ and $L^rv=(Lv)\sep$,
\par
$vL_0^r$ is the $p'$-divisible hull of $vL_0$ and $L_0^rv=(L_0v)\sep$.
\sn
It follows that
\begin{equation}                            \label{vLvL0}
(vL^r:vL_0^r)\>=\>(vL:vL_0)\;\;\;\mbox{ and }\;\;\; [L^rv:L_0^rv]\>=\>
[Lv:L_0v]\;.
\end{equation}
From (\ref{d=d}), (\ref{gr0}) and (\ref{vLvL0}), keeping in mind that
$L.K^r=L^r$ and $K^r=L_0^r$, we deduce
\begin{eqnarray*}
\mbox{\rm d}(L.K^r|K^r,v) & = & \mbox{\rm d}(L^r|L_0^r,v)\;=\;
\frac{[L^r:L_0^r]}{(vL^r:vL_0^r)[L^rv:L_0^rv]}\\
& = & \frac{[L:L_0]}{(vL:vL_0)\cdot [Lv:L_0v]}
\;=\; \mbox{\rm d}(L | L_0,v)\;=\;\mbox{\rm d}(L|K^h,v)\;.
\end{eqnarray*}
\end{proof}

\pars
We can now describe the ramification theoretic proof for the lemma of
Ostrowski (see also \cite[Corollary to Theorem~25, p.~78]{[Z--S]}). Take
a finite extension $(L'|L,v)$ of henselian fields. Then $L=L^h$. By the
foregoing theorem, ${\rm d}(L'|L,v)={\rm d}(L'.L^r|L^r, v)$. It follows
from Theorem~\ref{V} that $[L'.L^r:L^r]$ is a power of $p$. Hence also
${\rm d}(L'.L^r|L^r,v)$, being a divisor of it, is a power of $p$.

We see that non-trivial defects can only appear between $K^r$ and
$\tilde{K}$, or equivalently, between $K^h$ and $W$. These are the areas
of \bfind{wild ramification}, whereas the extension from $K^i$ to $K^r$
is the area of \bfind{tame ramification}. Hence, local uniformization in
characteristic 0 and the classification problem for valued fields of
residue characteristic 0 only have to deal with tame ramification, while
the two problems we described in Section~\ref{sect2p} also have to fight
the wild ramification.

\parm
An algebraic extension $(L'|L,v)$ of henselian fields is called
\textbf{unramified} if every finite subextension is unramified.
An algebraic extension $(L'|L,v)$ of henselian fields is called
\bfind{tame} if every finite subextension $(L''|L,v)$ is defectless and
such that $(L''v|Lv)$ is separable and $p$ does not divide $(vL'':vL)$.
A henselian field $(L,v)$ is called a \bfind{tame field} if $(\tilde{L}|
L, v)$ is a tame extension, and it is called a \textbf{separably tame
field} if $(L\sep|L,v)$ is a tame extension. The fields $W$ of
Theorem~\ref{pank} are examples of tame fields.

The proof of the following theorem is given in \cite{[Ku2]}.

\begin{theorem}                             \label{mue}
The absolute inertia field is the unique maximal unramified extension of
$K^h$ in $(\tilde{K},v)$. The absolute ramification field is the unique
maximal tame extension of $K^h$ in $(\tilde{K},v)$.
\end{theorem}

Note that an extension is tame if and only if it is defectless and
``tamely ramified'' in the sense of \cite[\$22]{[End]}. As we have
already mentioned, our notion of ``unramified'' is the same as
``defectless'' plus ``unramified'' in the sense of \cite[\$22]{[End]}
Hence for defectless valuations, the above theorem follows from
\cite[Corollary~(22.9)]{[End]}.

\parm
We summarize our main results in the following table:\n
\setlength{\unitlength}{0.0015\textwidth}
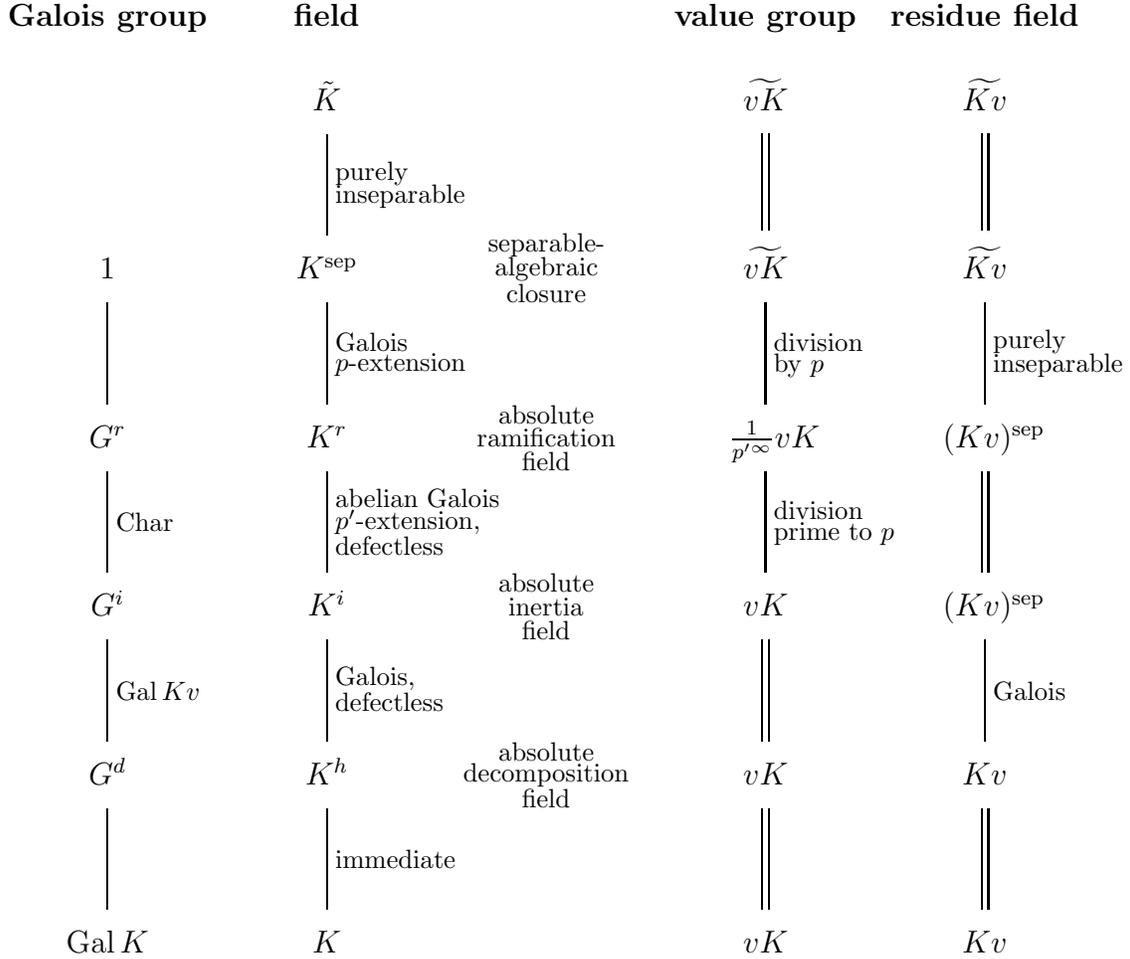
\begin{picture}(650,640)(10,0)              \label{ramtable}
\put(70,50){\bbox{$\Gal K$}}
\put(70,150){\bbox{$G\dec$}}
\put(70,250){\bbox{$G^i$}}
\put(70,350){\bbox{$G^r$}}
\put(70,450){\bbox{$1$}}
\put(70,550){\bbox{$$}}
\put(70,600){\bbox{\bf Galois group}}
\put(200,50){\bbox{$K$}}
\put(200,150){\bbox{$K^h$}}
\put(200,250){\bbox{$K^i$}}
\put(200,350){\bbox{$K^r$}}
\put(200,450){\bbox{$K\sep$}}
\put(200,550){\bbox{$\tilde{K}$}}
\put(200,600){\bbox{\bf {field}}}
\put(460,50){\bbox{$vK$}}
\put(460,150){\bbox{$vK$}}
\put(460,250){\bbox{$vK$}}
\put(460,350){\bbox{$\;\;\frac{1}{p'{}^\infty}vK$}}
\put(460,450){\bbox{$\widetilde{vK}$}}
\put(460,550){\bbox{$\widetilde{vK}$}}
\put(460,600){\bbox{\bf {value group}}}
\put(590,50){\bbox{$Kv$}}
\put(590,150){\bbox{$Kv$}}
\put(590,250){\bbox{$\;\;(Kv)\sep$}}
\put(590,350){\bbox{$\;\;(Kv)\sep$}}
\put(590,450){\bbox{$\widetilde{Kv}$}}
\put(590,550){\bbox{$\widetilde{Kv}$}}
\put(590,600){\bbox{\bf {residue field}}}
\put(70,130){\line(0,-1){60}}    
\put(70,230){\line(0,-1){60}}    
\put(70,330){\line(0,-1){60}}    
\put(70,430){\line(0,-1){60}}    
\put(200,130){\line(0,-1){60}}    
\put(200,230){\line(0,-1){60}}    
\put(200,330){\line(0,-1){60}}    
\put(200,430){\line(0,-1){60}}    
\put(200,530){\line(0,-1){60}}    
\put(458,130){\line(0,-1){60}}    
\put(462,130){\line(0,-1){60}}    
\put(458,230){\line(0,-1){60}}    
\put(462,230){\line(0,-1){60}}    
\put(460,330){\line(0,-1){60}}    
\put(460,430){\line(0,-1){60}}    
\put(458,530){\line(0,-1){57}}    
\put(462,530){\line(0,-1){57}}    
\put(588,130){\line(0,-1){60}}    
\put(592,130){\line(0,-1){60}}    
\put(590,230){\line(0,-1){60}}    
\put(588,330){\line(0,-1){60}}    
\put(592,330){\line(0,-1){60}}    
\put(590,430){\line(0,-1){60}}    
\put(588,530){\line(0,-1){57}}    
\put(592,530){\line(0,-1){57}}    
\put(75,200){\makebox(0,0)[l]{\footnotesize\rm $\Gal Kv$}}%
\put(75,300){\makebox(0,0)[l]{\footnotesize\rm Char}}%
\put(205,100){\makebox(0,0)[l]{\footnotesize\rm immediate}}%
\put(205,207){\makebox(0,0)[l]{\footnotesize\rm Galois,}}%
\put(205,193){\makebox(0,0)[l]{\footnotesize\rm defectless}}%
\put(205,314){\makebox(0,0)[l]{\footnotesize\rm abelian Galois}}%
\put(205,300){\makebox(0,0)[l]{\footnotesize\rm $p'$-extension,}}%
\put(205,286){\makebox(0,0)[l]{\footnotesize\rm defectless}}%
\put(205,407){\makebox(0,0)[l]{\footnotesize\rm Galois}}%
\put(205,393){\makebox(0,0)[l]{\footnotesize\rm $p$-extension}}%
\put(205,507){\makebox(0,0)[l]{\footnotesize\rm purely}}%
\put(205,493){\makebox(0,0)[l]{\footnotesize\rm inseparable}}%
\put(465,407){\makebox(0,0)[l]{\footnotesize\rm division}}%
\put(465,393){\makebox(0,0)[l]{\footnotesize\rm by $p$}}%
\put(465,307){\makebox(0,0)[l]{\footnotesize\rm division}}%
\put(465,293){\makebox(0,0)[l]{\footnotesize\rm prime to $p$}}%
\put(595,200){\makebox(0,0)[l]{\footnotesize\rm Galois}}%
\put(595,407){\makebox(0,0)[l]{\footnotesize\rm purely}}%
\put(595,393){\makebox(0,0)[l]{\footnotesize\rm inseparable}}%
\put(330,164){\makebox(0,0)[c]{\footnotesize\rm absolute}}%
\put(330,150){\makebox(0,0)[c]{\footnotesize\rm decomposition}}%
\put(330,136){\makebox(0,0)[c]{\footnotesize\rm field}}%
\put(330,264){\makebox(0,0)[c]{\footnotesize\rm absolute}}%
\put(330,250){\makebox(0,0)[c]{\footnotesize\rm inertia}}%
\put(330,236){\makebox(0,0)[c]{\footnotesize\rm field}}%
\put(330,364){\makebox(0,0)[c]{\footnotesize\rm absolute}}%
\put(330,350){\makebox(0,0)[c]{\footnotesize\rm ramification}}%
\put(330,336){\makebox(0,0)[c]{\footnotesize\rm field}}%
\put(330,464){\makebox(0,0)[c]{\footnotesize\rm separable-}}%
\put(330,450){\makebox(0,0)[c]{\footnotesize\rm algebraic}}%
\put(330,436){\makebox(0,0)[c]{\footnotesize\rm closure}}%
\end{picture}
\n
where $\frac{1}{p'{}^\infty}vK$ denotes the $p'$-divisible hull of
$vK$ and Char denotes the character group (\ref{cg}).

\pars
In algebraic geometry, the absolute inertia field is often
called the \bfind{strict henselization}. Theorem~\ref{lu=ig} can be
understood as saying that the Implicit Function Theorem, or
equivalently, Hensel's Lemma, works within and only within the strict
henselization. That the limit is the strict henselization and not the
henselization becomes intuitively clear when one considers one of the
equivalent forms of Hensel's Lemma which states that if $f$ has
coefficients in the valuation ring of a henselian field, then every
simple root of the reduced polynomial $fv$ (obtained by replacing the
coefficients by their residues) can be lifted to a root of $f$. On the
other hand, irreducible polynomials have only simple roots if and only
if they are separable. Hence it is clear that Hensel's Lemma works as
long as the residue field extensions are separable, which is the case
between $K^h$ and $K^i$.

\end{section}

%
%
\begin{section}{Two theorems}
%
%
%
\begin{subsection}{The Stability Theorem}   \label{subsst}
In this section we present two theorems about the defect which we have
used for our results on local uniformization and in the model theory of
valued fields in positive characteristic. The first one describes
situations where no defect appears. The second one deals with with
certain situation where defect may well appear, but shows that the
defect can be eliminated.

Let $(L|K,v)$ be an extension of valued fields of finite transcendence
degree. Then the following well known form of the ``Abhyankar
inequality'' holds:
\begin{equation}                            \label{wtdgeq}
\trdeg L|K \>\geq\> \rr vL/vK \,+\, \trdeg Lv|Kv\;,
\end{equation}
where $\rr vL/vK:=\dim_{\Q}\,(vL/vK)\otimes\Q$ is the \bfind{rational
rank} of the abelian group $vL/vK$, i.e., the maximal number of
rationally independent elements in $vL/vK$. This inequality is a
consequence of Theorem~1 of \cite[Chapter VI, \S10.3]{[B]}, which
states that if
\begin{equation}                            \label{vtb}
\left\{
\begin{array}{l}
x_1,\ldots,x_\rho,y_1,\ldots,y_\tau\in L\mbox{ such that}\\
vx_1,\ldots,vx_\rho \mbox{ are
rationally independent over $vK$, and}\\
y_1v,\ldots,y_\tau v \mbox{ are algebraically independent over $Kv$,}
\end{array}\right.
\end{equation}
then $x_1,\ldots,x_\rho,y_1,\ldots,y_\tau$ are algebraically independent
over $K$. We will say that $(L|K,v)$ is \textbf{without transcendence
defect}\index{transcendence defect} if equality holds in (\ref{wtdgeq}).
In this case, every set $\{x_1,\ldots,x_\rho,y_1,\ldots,y_\tau\}$
satisfying (\ref{vtb}) with $\rho=\rr vL/vK$ and $\tau=\trdeg Lv|Kv$ is
a transcendence basis of $L|K$.

If $(F|K,v)$ is a valued function field without transcendence defect,
then the extensions $vF|vK$ and $Fv|Kv$ are finitely generated (cf.\
\cite[Corollary~2.2]{[Kn-Ku1]}).

\begin{theorem}                                    \label{ai}
{\bf (Generalized Stability Theorem)}\index{Stability Theorem}\n
Let $(F|K,v)$ be a valued function field without transcendence
defect. If $(K,v)$ is a defectless field, then $(F,v)$ is a defectless
field. The same holds for ``inseparably defectless'' in the place of
``defectless''. If $vK$ is cofinal in $vF$, then it also holds for
``separably defectless'' in the place of ``defectless''.
\end{theorem}

If the base field $K$ is not a defectless field, we can say at least the
following:
\begin{corollary}                                \label{aife}
Let $(F|K,v)$ be a valued function field without transcendence
defect, and $E|F$ a finite extension. Fix an extension of $v$ from $F$
to $\tilde{K}.F\,$. Then there is a finite extension $L_0|K$ such that
for every algebraic extension $L$ of $K$ containing $L_0\,$, $(L.F,v)$
is defectless in $L.E\,$. If $(K,v)$ is henselian, then $L_0|K$ can be
chosen to be purely wild.
\end{corollary}

Theorem~\ref{ai} was stated and proved in \cite{[Ku1]}; the proof
presented in \cite{[Ku8]} is an improved version.

The theorem has a long and interesting history. Hans Grauert and
Reinhold Remmert \cite[p.~119]{[Gra--Re]} first proved it in a very
restricted case, where $(K,v)$ is an algebraically closed complete
discretely valued field and $(F,v)$ is discrete too. A generalization
was given by Laurent Gruson \cite[Th\'eor\`eme~3, p.~66]{[Gru]}. A good
presentation of it can be found in the book on non-archi\-me\-dean
analysis by Siegfried Bosch, Ulrich G\"untzer and Reinhold Remmert
\cite[\S 5.3.2, Theorem 1]{[Bos--Gu--Re]}. Further generalizations are
due to Michel Matignon and Jack Ohm, and also follow from results in
\cite{[Gre--M--Pop1]} and \cite{[Gre--M--Pop2]}. Ohm arrived
independently of \cite{[Ku1]} at a general version of the Stability
Theorem for the case of $\trdeg L|K =\trdeg Lv|Kv$ (see the second
theorem on p.~306 of \cite{[Oh]}). He deduces his theorem from
Proposition 3 on p.\ 215 of \cite{[Bos--Gu--Re]}, (more precisely, from
a generalized version of this proposition which is proved but not stated
in \cite{[Bos--Gu--Re]}).

All authors mentioned in the last paragraph use methods of
non-archime\-dean analysis, and all results are restricted to the case
of $\trdeg F|K=\trdeg Fv|Kv$. In this case we call the extension
$(F|K,v)$ \textbf{residually transcendental}, and we call the valuation
$v$ a \bfind{constant reduction} of the algebraic function field $F|K$.
The classical origin of such valuations is the study of curves over
number fields and the idea to reduce them modulo a $p$-adic valuation.
Certainly, the reduction should again render a curve, this time over a
finite field. This is guaranteed by the condition $\trdeg F|K=\trdeg
Fv|Kv$, where $F$ is the function field of the curve and $Fv$ will be
the function field of its reduction. Naturally, one seeks to relate the
genus of $F|K$ to that of $Fv|Kv\,$. Several authors proved \textbf{genus
inequalities} (see, for example, \cite{[De],[M],[Gre--M--Pop1]} and the
survey given in \cite{[Gre]}). To illustrate the use of the defect, we
will cite an inequality proved by Barry Green, Michel Matignon and
Florian Pop \cite[Theorem~3.1]{[Gre--M--Pop1]}. Let $F|K$ be a function
field of transcendence degree 1, and $v$ a constant reduction of $F|K$.
We choose a henselization $F^h$ of $(F,v)$; all henselizations of
subfields of $F$ will be taken in $F^h$. We wish to define a defect of
the extension $(F^h|K^h,v)$ even though this extension is not algebraic.
The following result helps:

\begin{theorem}                             \label{IT}
{\bf (Independence Theorem)}\n
The defect of the algebraic extension $(F^h|K(t)^h,v)$ is independent of
the choice of the element $t\in F$, provided that $tv$ is transcendental
over $Kv$.
\end{theorem}
In \cite{[Oh]}, Ohm proves a more general version of this theorem for
arbitrary transcendence degree, using his version of the Stability
Theorem. The Stability Theorem tells us that in essence, the defect of a
residually transcendental function field, and more generally, of a
function field without transcendence defect, can only come from the base
field. The following general Independence Theorem was proved in
\cite[Theorem~5.4 and Corollary~5.6]{[Ku1]}:

\begin{theorem}
Take a valued function field $(F|K,v)$ without transcendence defect, and
set $\rho=\rr vF/vK$ and $\tau=\trdeg Fv|Kv$. The defect of the
extension $(F^h|K(x_1, \ldots,x_\rho,y_1,\ldots,y_\tau)^h,v)$ is
independent of the choice of the elements\linebreak
$x_1,\ldots,x_\rho,y_1,
\ldots,y_\tau$ as long as they satisfy (\ref{vtb}). Moreover, there is a
finite extension $K'|K$ such that $(F^h.K'|K'(x_1,
\ldots,x_\rho,y_1,\ldots,y_\tau)^h,v)$ is defectless and
\[
{\rm d}(F^h.K'|K(x_1,\ldots,x_\rho,y_1,\ldots,y_\tau)^h,v)\>=\>
{\rm d}(K^h.K'|K^h,v)\;.
\]
\end{theorem}
A special case for simple transcendental extensions $(K(x)|K,v)$
satisfying the condition $\trdeg K(x)|K = \rr vK(x)/vK$ was proved by
Sudesh Khanduja in \cite{[KH10]}.

\pars
An interesting proof of Theorem~\ref{IT} is given in
\cite{[Gre--M--Pop1]}, as it introduces another notion of defect. We
take any valued field extension $(L|K,v)$ and a finite-dimensional
$K$-vector space $V\subseteq L$. We choose a system $\mathcal V$ of
representatives of the cosets $va+vK$, $0\ne a\in V$. For every
$K$-vector space $W\subseteq V$ and every $\gamma\in\mathcal V$ we set
$W_\gamma:=\{a\in W\mid va\geq 0\}$ and $W_\gamma^\circ:=\{a\in W\mid
va>0\}$. The quotient $W_\gamma/W_\gamma^\circ$ is in a natural way a
$Kv$-vector space. The \bfind{vector space defect} of $(V|K,v)$ is
defined as
\[
{\rm d}^{\rm vs}(V|K,v)\>:=\> \sup_{W\subseteq V}
\frac{\dim_K W}{\sum_{\gamma\in\mathcal V}\dim_{Kv}
W_\gamma/W_\gamma^\circ}\;,
\]
where the supremum runs over all finite-dimensional subspaces $W$. For a
finite extension $(L|K,v)$, by \cite[Proposition~2.2]{[Gre--M--Pop1]},
\[
{\rm d}^{\rm vs}(L|K,v)\>=\> \frac{[L:K]}{(vL:vK)[Lv:Kv]}
\]
which is equal to the ordinary defect ${\rm d}(L|K,v)$ if the extension
of $v$ from $K$ to $L$ is unique.

Note that quotients of the form $W_\gamma/W_\gamma^\circ$ also appear in
the definition of the graded ring of a subring in a valued field, then
often written as ``${\mathcal P}_\gamma/{\mathcal P}_\gamma^+$'' (see,
for instance, \cite[\$2]{[V2]}). Graded rings are used by Bernard
Teissier in his program for a ``characteristic blind'' local
uniformization, see \cite{[Te]}.

\pars
The following result (\cite[Theorem~2.13]{[Gre--M--Pop1]}) implies
the Independence Theorem~\ref{IT}:

\begin{theorem}                             \label{IT1}
For every element $t\in F$ such that $tv$ is transcendental over
$Kv$,
\[
{\rm d}^{\rm vs}(F|K,v)\>=\> (F^h|K(t)^h,v)\;.
\]
\end{theorem}

Now we are ready to cite the genus inequality for an algebraic function
field $F|K$ with distinct constant reductions $v_1,\ldots,v_s$ which
have a common restriction to $K$. We assume in addition that $K$
coincides with the constant field of $F|K$ (the relative algebraic
closure of $K$ in $F$). Then:
\begin{equation}
1-g_F\leq 1-s+\sum_{i=1}^{s} {\rm d}_i {\rm e}_i r_i (1-g_i)
\end{equation}
where $g_F$ is the genus of $F|K$ and $g_i$ the genus of $Fv_i|Kv_i\,$,
$r_i$ is the degree of the constant field of $Fv_i|Kv_i$ over $Kv_i\,$,
${\rm d}_i={\rm d}^{\rm vs}(F|K,v_i)$, and ${\rm e}_i= (v_iF:v_iK)$
(which is always finite in the constant reduction case, see, for
instance, \cite[Corollary~2.7]{[Ku6]}). It follows that constant
reductions $v_1,v_2$ with common restriction to $K$ and $g_1=g_2=
g_F\geq 1$ must be equal. In other words, for a fixed valuation on $K$
there is at most one extension $v$ to $F$ which is a \textbf{good
reduction}, that is, (i) $g_F=g_{Fv}\,$, (ii) there exists $f\in F$ such
that $vf=0$ and $[F:K(f)]= [Fv:Kv(fv)]$, (iii) $Kv$ is the constant
field of $Fv|Kv\,$. An element $f$ as in (ii) is called a
\textbf{regular function}.

More generally, $f$ is said to have the \textbf{uniqueness property} if
$fv$ is transcendental over $Kv$ and the restriction of $v$ to $K(f)$
has a unique extension to $F$. In this case, $[F:K(f)]={\rm d}\cdot
{\rm e}\cdot [Fv:Kv(fv)]$ where d is the defect of $(F^h|K^h,v)$ and
e$\>=(vF:vK(f))=(vF:vK)$. If $K$ is algebraically closed, then e$\>=1$,
and it follows from the Stability Theorem that d$\>=1$; hence in this
case, every element with the uniqueness property is regular.

It was proved in \cite[Theorem~3.1]{[Gre--M--Pop2]} that $F$ has an
element with the uniqueness property already if the restriction of $v$
to $K$ is henselian. The proof uses Abraham Robinson's model
completeness result for algebraically closed valued fields, and
ultraproducts of function fields. Elements with the uniqueness property
also exist if $vF$ is a subgroup of $\Q$ and $Kv$ is algebraic over a
finite field. This follows from work in \cite{[Gre--M--Pop3]} where the
uniqueness property is related to the \textbf{local Skolem property}
which gives a criterion for the existence of algebraic $v$-adic integral
solutions on geometrically integral varieties. This result is a special
case of a theorem proved in \cite{[KH6]} which states that elements with
the uniqueness property exist if and only if the completion of $(K,v)$
is henselian.

As an application to rigid analytic spaces, the Stability Theorem is
used to prove that the quotient field of the free Tate algebra $T_n(K)$
is a defectless field, provided that $K$ is. This in turn is used to
deduce the \textbf{Grauert--Remmert Finiteness Theorem}, in a generalized
version due to Gruson; see \cite[pp.~214--220]{[Bos--Gu--Re]} for ``a
simplified version of Gruson's approach''.

\pars
In contrast to the approaches that use methods of non-archimedean
ana\-lysis, we give in \cite{[Ku1]}, \cite{[Ku2]} and \cite{[Ku8]} a new
proof which replaces the analytic methods by valuation theoretical
arguments. Such arguments seem to be more adequate for a theorem that is
of (Krull) valuation theoretical nature.

Our approach has much in common with Abhyankar's method of using
ramification theory in order to reduce the question of resolution of
singularities to the study of purely inseparable extensions and of
Galois extensions of degree $p$ and the search for suitable normal forms
of Artin-Schreier-like minimal polynomials (cf.\ \cite{[Ab2]}). Given a
finite separable extension $(L'|L,v)$ of henselian fields of positive
characteristic, we can study its properties by lifting it up to the
absolute ramification fields. From Lemma~\ref{dram} we know that the
defect of $(L'|L,v)$ is equal to the defect of $(L'.L^r|L^r,v)$. From
Lemma~\ref{RemASe} we know that the extension $L'.L^r|L^r$ is a tower of
Artin-Schreier extensions.

Abhyankar's ramification thoretical reduction to Artin-Schreier
extensions and purely inseparable extensions is also used by Vincent
Cossart and Olivier Piltant in \cite{[CoP1]} to reduce resolution of
singularities of threefolds in positive characteristic to local
uniformization on Artin-Schreier and purely inseparable coverings. The
Artin-Schreier extensions appearing through this reduction are not
necessarily defect extensions. According to Piltant, those that are, are
harder to treat than the defectless ones.

In the situation
of Theorem~\ref{ai}, we have to prove for $L=F^h$ that $(L'|L,v)$ is
defectless, or equivalently, that each Artin-Schreier extension in the
tower is defectless. Looking at the first one in the tower, assume that
it is generated by a root $\vartheta$ of a polynomial $X^p-X-a$ with
$a\in L^r$. Using the additivity of the Frobenius in characteristic $p$,
we see that the element $\vartheta-c$, which generates the same
extension, has minimal polynomial $X^p-X-(a-c^p+c)$. Hence if $a$
contains some $p$-th power $c^p$, we can replace it by $c$ without
changing the extension. Using this fact and the special structure of
$(F^h,v)$ given by the assumptions of Theorem~\ref{ai} on $(F,v)$, we
deduce normal forms for $a$ which allow us to read off that the
extension is defectless. This fact in turn implies that
$(F^h(\vartheta),v)$ is again of the same special form as $(F^h,v)$,
which enables us to proceed by induction over the extensions in the
tower.

Note that when algebraic geometers work with Artin-Schreier extensions
they usually work with polynomials of the form $X^p-dX-a$. The reason is
that they work over rings and not over fields. A polynomial like
$X^p-b^{p-1}X-a$ over a ring $R$ can be transformed into the polynomial
$X^p-X-a/b^p$, as we have seen in Example~\ref{transf}, but $a/b^p$ does
in general not lie in the ring anymore. Working with a polynomial of the
form $X^p-X-a$ is somewhat easier than with a polynomial of the form
$X^p-dX-a$, and it suffices to derive normal forms as needed for the
proof of Theorem~\ref{ai}, and of Theorem~\ref{hr} which we will discuss
below.

In the case of mixed characteristic, where the valued fields have
characteristic 0 and their residue fields have positive characteristic,
Artin-Schreier extensions are replaced by Kummer extensions (although
re-written with corresponding Artin-Schreier polynomials), and
additivity is replaced by quasi-additivity (cf.\ Lemma~\ref{quasadd}).

Related normal form results can be found in the work of Helmut Hasse,
George Whaples, and in Matignon's proof of his version of
Theorem~\ref{ai}. See also Helmut Epp's paper \cite{[Ep]}, in particular
the proof of Theorem~1.3. This proof contains a gap which was filled in
\cite{[KuE]}.

\parm
Let us reconsider Examples~\ref{ex1} and~\ref{ex2} in the light of
Theorem~\ref{ai}. In Example~\ref{ex1} we have an extension without
transcendence defect if and only the transcendence degree is 1. In this
case, $(K(t)^h,v)$ is a defectless field, and we have that $F\subset
K(t)^h$. In the case of higher transcendence degree, this may not be the
case, as Example~\ref{exampFKS} shows. At least we know that every
separable extension of $K(T)^h$ is defectless since it is discretely
valued. The situation is different in Example~\ref{ex2}. If the
extension is without transcendence defect, then again, $(K(t)^h,v)$ is
a defectless field, and moreover, $vF/vK$ and $Fv|Kv$ are finitely
generated (\cite[Corollary~2.7]{[Ku6]}). But if $\chara K>0$, then there
are valuations $v$ on $K(x,y)$, trivial on $K$, such that $K(x,y)v=K$,
$vK(x,y)$ not finitely generated, and such that $(K(x,y),v)$ admits an
infinite tower of Artin-Schreier defect extensions
(\cite[Theorem~1.2]{[Ku6]}).

\parm
Applications of Theorem~\ref{ai} are:
\sn
$\bullet$ \ {\bf Elimination of ramification.} \
In \cite{[Kn-Ku1]} we use Theorem~\ref{ai} to prove:

\begin{theorem}                             \label{hrwtd}
Take a defectless field $(K,v)$ and a valued function field $(F|K,v)$
without transcendence defect. Assume that $Fv|Kv$ is a separable
extension and $vF/vK$ is torsion free. Then $(F|K,v)$ admits elimination
of ramification in the following sense: there is a transcendence basis
$T=\{x_1,\ldots,x_r,y_1,\ldots, y_s\}$ of $(F|K,v)$ such that
\pars
a) \ $vF= vK\oplus\Z vx_1\oplus \ldots\oplus
\Z vx_r$,\par
b) \ $y_1 v,\ldots,y_s v$ form a separating transcendence
basis of $Fv|Kv$.
\sn
For each such transcendence basis $T$ and every extension of $v$
to the algebraic closure of $F$, $(F^h|K(T)^h,v)$ is unramified.
\end{theorem}

\begin{corollary}                             \label{elram}
Let $(F|K,v)$ be a valued function field without transcendence
defect. Fix an extension of $v$ to $\tilde{F}$. Then there is a finite
extension $L_0|K$ and a transcendence basis $T$ of $(L_0.F|L_0,v)$ such
that for every algebraic extension $L$ of $K$ containing $L_0\,$, the
extension $((L.F)^h|L(T)^h,v)$ is unramified.
\end{corollary}

\mn
$\bullet$ \ {\bf Local uniformization in positive and in mixed
characteristic.} \ We consider places $P$ and their associated
valuations $v=v_P$ of a function field $F|K$, by which we mean that
$P|_K$ is the identity and hence $v|_K$ is trivial. We write $aP=av$
and denote by ${\cal O}$ the valuation ring of $v$ on $F$. Rewriting
our earlier definition, we say that
$P$ admits \bfind{smooth local uniformization} if there is a model for
$F$ on which $P$ is centered at a smooth point, that is, if there are
$x_1,\ldots,x_n\in {\cal O}$ such that $F=K(x_1,\ldots,x_n)$ and the
point $x_1P,\ldots,x_nP$ is smooth. (Note that in \cite{[Kn-Ku1]} and
\cite{[Kn-Ku2]} we add a further condition, which we drop here for
simplicity.) The place $P$ is called an \bfind{Abhyankar place} if
equality holds in the Abhyankar inequality, which in the present case
means that $\trdeg F|K=\rr vF+\trdeg FP|K$.

Theorem~\ref{hrwtd} is a crucial ingredient for the following result
(cf.\ \cite[Theorem~1.1]{[Kn-Ku1]}, \cite{[Ku3]}):

\begin{theorem}                             \label{lu1}
Assume that $P$ is an Abhyankar place of the function field $F|K$ such
that $FP|K$ is separable. Then $P$ admits smooth local uniformization.
\end{theorem}
The analogous arithmetic case (\cite[Theorem~1.2]{[Kn-Ku1]}) uses
Theorem~\ref{ai} in mixed characteristic. Note that the condition
``$FP|K$ is separable'' is necessary since it is implied by elimination
of ramification.

\bn
$\bullet$ \ {\bf Model theory of valued fields.} \
In \cite{[Ku10]} we use Theorem~\ref{hrwtd} to prove the following
Ax--Kochen--Ershov Principle:
\begin{theorem}                             \label{AKEwtd}
Take a henselian defectless valued field $(K,v)$ and an extension
$(L|K,v)$ of finite transcendence degree without transcendence defect.
If $vK$ is existentially closed in $vL$ and $Kv$ is existentially
closed in $Lv$, then $(K,v)$ is existentially closed in $(L,v)$.
\end{theorem}
Let us continue our discussion from the end of Section~\ref{subsee}. The
conditions ``$vK\ec vL$'' and ``$Kv\ec Lv$'' imply that $vF/vK$ is
torsion free and $Fv|Kv$ is a separable extension. So we can apply
Theorem~\ref{hrwtd} to obtain the transcendence basis
$T=\{x_1,\ldots,x_r,y_1,\ldots, y_s\}$ of $(F|K,v)$ with the properties
as specified in that theorem. Because of these properties, the
embeddings of $vL$ in $v^*K^*$ and of $Lv$ in $K^*v^*$ lift to an
embedding $\iota_0$ of $(K(T),v)$ in $(K^*,v^*)$ over $K$. Using
Hensel's Lemma and the fact that $Fv|K(T)v$ is separable, one finds in
$F^h|K(T)$ a subextension $F_1|K(T)$ with $F_1v=Fv$ and $[F_1:K(T)]=
[F_1v:K(T)v]$ ($Fv|K(T)v$ is finite by the remark preceding
Theorem~\ref{ai}). Using Hensel's Lemma and the embedding of $Lv$ in
$K^*v^*$ again, one extends $\iota_0$ to an embedding $\iota_1$ of
$(F_1,v)$ in $(K^*,v^*)$. The extension $(F^h|F_1^h,v)$ is immediate,
and as it is an extension inside the unramified extension
$(F^h|K(T)^h,v)$, it must be defectless and hence trivial. As
$(K^*,v^*)$ is henselian, being an elementary extension of the henselian
field $(K,v)$, one can now use the universal property of henselizations
to extend $\iota_1$ to an embedding $\iota_2$ of $(F^h,v)$ in
$(K^*,v^*)$. The restriction of $\iota_2$ to $F$ is the desired
embedding which transfers every existential sentence valid in $(F,v)$
into $(K^*,v^*)$.

\end{subsection}

%
%
\begin{subsection}{Henselian Rationality of Immediate Function Fields}
Let us return to Example~\ref{exclhr}. If $(F,v)$ does not lie in the
henselization $K(x)^h$, we are lost. This happens if and only if
$(F^h|K(x)^h,v)$ has non-trivial defect (the equivalence holds because
$(F|K(x),v)$ is finite and immediate, $F^h=F.K(x)^h$ and henselizations
are immediate extensions).

So the question arises: how can we avoid the defect in the case of
immediate extensions? The answer is a theorem proved in \cite{[Ku1]}
(cf.\ \cite{[Ku2]} and \cite{[Ku11]}):
\begin{theorem}                \label{hr}
{\bf (Henselian Rationality)} \n
Let $(K,v)$ be a tame field and $(F|K,v)$ an immediate function field
of trans\-cendence degree 1. Then
\begin{equation}                            \label{hre}
\mbox{there is $x\in F$ such that }\; (F^h,v)\,=\,(K(x)^h,v)\;,
\end{equation}
that is, $(F|K,v)$ is henselian generated. The same holds over a
separably tame field $(K,v)$ if in addition $F|K$ is separable.
\end{theorem}
Since the assertion says that $F^h$ is equal to the henselization of a
rational function field, we also call $F$ \bfind{henselian rational} in
this case.
For valued fields of residue characteristic 0, the assertion is a direct
consequence of the fact that every such field is defectless. Indeed,
take any $x\in F\setminus K$. Then $K(x)|K$ cannot be algebraic since
otherwise, $(K(x)|K,v)$ would be a proper finite immediate (and hence
defect) extension of the tame field $(K,v)$, a contradiction to the
definition of ``tame''. Hence, $F|K(x)$ is algebraic and immediate.
Therefore, $(F^h|K(x)^h,v)$ is algebraic and immediate too. But since it
cannot have a non-trivial defect, it must be trivial. This proves that
$(F,v)\subset (K(x)^h,v)$. In contrast to this, in the case of positive
residue characteristic only a very carefully chosen $x\in F\setminus K$
will do the job. As for the Generalized Stability Theorem, the proof of
Theorem~\ref{hr} in positive characteristic uses ramification theory and
the deduction of normal forms for Artin-Schreier extensions. This time
however, all Artin-Schreier extensions are immediate and hence defect
extensions. The normal forms serve a different purpose, namely, to
find a suitable generator $x$. The proof also uses significantly
a theory of immediate extensions which builds on Kaplansky's paper
\cite[Sections~2 and~3]{[Ka]}.

\sn
{\bf Open problem (HR):} Improve Theorem~\ref{hr} by finding versions
that work with weaker assumptions. For instance, can the assumption
``tame'' be replaced by ``henselian and perfect'' or just ``perfect'',
or can it even be dropped altogether? Then, even with a weaker
assumption on $(K,v)$, can the assumption ``immediate'' be replaced by
``$vF/vK$ is a torsion group and $Fv=Kv$''?

Note that in order to allow $Fv|Kv$ to be any algebraic extension, a
possible generalization of Theorem~\ref{hr} would have to replace
(\ref{hre}) by
\begin{equation}
\mbox{there is $x\in F$ such that }\; (F^i,v)\,=\,(K(x)^i,v)\;.
\end{equation}

\parm
Applications of Theorem~\ref{hr} in conjunction with Theorem~\ref{ai}
are:
\sn
$\bullet$ \ {\bf Local uniformization in positive and in mixed
characteristic.} \ Theorem~\ref{hr} together with Theorem~\ref{hrwtd} is
a crucial ingredient for the proof of ``local uniformization by
alteration'' (cf.\ \cite[Theorem~1.2]{[Kn-Ku2]}, \cite{[Ku3]}):

\begin{theorem}                             \label{lu2}
Assume that $P$ is a place of the function field $F|K$. Then there is a
finite extension $F'|F$ and an extension $P'$ of $P$ from $F$ to $F'$
such that $P'$ admits smooth local uniformization. The extension $F'|F$
can be chosen to be Galois. Alternatively, it can be chosen such that
$(F',P')|(F,P)$ is purely wild, hence $v_{P'}F'/v_P F$ is a $p$-group
and $F'P'|FP$ is purely inseparable.
\end{theorem}
The analogous arithmetic case (\cite[Theorem~1.4]{[Kn-Ku2]})
uses Theorems~\ref{hr} and ~\ref{ai} in mixed characteristic.
While local uniformization by alteration follows from de Jong's
resolution of singularities by alteration (see \cite{[AdJ]}), the
additional information on the extension $F'|F$ does not follow.
Moreover, the proofs of Theorems~\ref{lu1} and~\ref{lu2} use only
valuation theory.

Recently, Michael Temkin (\cite[Corollary~1.3]{[T]}) proved
``Inseparable Local Uniformization'':

\begin{theorem}                             \label{lu3}
In the setting of Theorem~\ref{lu2}, the extension $F'|F$ can also be
chosen to be purely inseparable.
\end{theorem}

It is interesting that local uniformization has now been proved up to
separable alteration on the one hand, and up to purely inseparable
alteration on the other. These two results are somewhat ``orthogonal''
to each other. Can they be put together to get rid of alteration? While
this appears to be an attractive thought at first sight, one should keep
in mind Example~\ref{transf} which shows that every purely inseparable
defect extension of degree $p$ of $(L,v)$ which does not lie in the
completion of $(L,v)$ can be transformed into an Artin-Schreier defect
extension. Thus, the ``same'' defect may appear in a separable extension
and in a purely inseparable extension (see the next section for
details), which leaves us the choice to kill it either with separable or
with inseparable alteration. So this fact does not in itself indicate
whether we need or do not need alteration for local uniformization.

\n
{\bf Open problem (LU):} Prove (or disprove) local uniformization
without extension of the function field.

In fact, one reason for the extension of the function field in our
approach is the fact that we apply Theorem~\ref{hr} to fields of lower
transcendence degree than the function field itself. However,
subfunction fields are too small to be tame fields, so we enlarge our
intermediate fields so that they become (separably) tame, and once we
have found
local uniformization in this larger configuration, we collect the only
finitely many new elements that are needed for it and adjoin them to the
original function field. So we see that if we can weaken the assumptions
of Theorem~\ref{hr}, then possibly we will need smaller extensions of
our function field. Temkin's work contains several developments in
this direction, one of which we will discuss in more detail in the next
section.

\bn
$\bullet$ \ {\bf Model theory of valued fields.} \
In \cite{[Ku10]} we use Theorem~\ref{hr} together with
Theorem~\ref{hrwtd} to prove the following:
\begin{theorem}                             \label{AKEtame}
\n
a) If $(K,v)$ is a tame field, then the Ax--Kochen--Ershov Principle
(\ref{AKE}) holds.
\sn
b) The Classification Problem for valued fields has a positive solution
for tame fields:
If $(K,v)$ and $(L,v)$ are tame fields such that $vK$ and $vL$ are
elementarily equivalent as fields and $Kv$ and $Lv$ are elementarily
equivalent as ordered groups, then $(K,v)$ and $(L,v)$ are elementarily
equivalent as valued fields.
\end{theorem}
\n
This theorem comprises several classes of valued fields for which the
classification had already been known to hold, such as the already
mentioned henselian fields with residue fields of characteristic 0.

\sn
{\bf Open problem (AKE):} Prove Ax--Kochen--Ershov Principles for
classes of non-perfect valued fields of positive characteristic. This
problem is connected with the open problem whether the elementary theory
of $\Fp((t))$ is decidable (cf.\ \cite{[Ku2]}, \cite{[Ku5]},
\cite{[Ku10]}).

\end{subsection}
\end{section}

%
\begin{section}{Two types of Artin-Schreier defect extensions}
In this section, we assume all fields to have characteristic $p>0$.
In Section~\ref{sectexadef} we have given several examples of
Artin-Schreier defect extensions, i.e., Artin-Schreier defect extensions
with non-trivial defect. Note that every such extension is immediate.
Some of our examples were derived from immediate purely inseparable
extensions of degree $p$ (Examples~\ref{transf} and~\ref{transfdt}). If
an Artin-Schreier defect extension is derived from a purely inseparable
defect extension of degree $p$ as in Example~\ref{transf}, then we call
it a \bfind{dependent Artin-Schreier defect extension}. If it cannot be
derived in this way, then we call it an \bfind{independent
Artin-Schreier defect extension}. More precisely, an Artin-Schreier
defect extension $(L(\vartheta)|L,v)$ with $\vartheta^p-\vartheta\in L$
is defined to be dependent if there is a purely inseparable extension
$(L(\eta)|L,v)$ of degree $p$ such that
\[
\mbox{for all }c\in L,\;\; v(\vartheta-c)\>=\> v(\eta-c)\;.
\]

The extension $(L(\vartheta)|L,v)$ constructed in Example~\ref{fse} is
an independent Artin-Schreier defect extension. This is obvious if we
choose $K=\Fp(t)$ or $K=\Fp((t))$ because then $L$ is the perfect hull
of $K$ and does not admit any purely inseparable defect extensions at
all. But if for instance, $K=\Fp(t,s)$ with $s\in\Fp((t))$
transcendental over $\Fp(t)$, then $L$ is not perfect. How do we know
then that $(L(\vartheta)|L,v)$ is an independent Artin-Schreier defect
extension? The answer is given by the following characterization proved
in \cite{[Ku1]} (see \cite{[Ku2]} and \cite{[Ku9]}):

\begin{theorem}                             \label{dindcr}
Take an Artin-Schreier defect extension $(L(\vartheta)|L,v)$ with
$\vartheta^p-\vartheta\in L$. Then this extension is independent if and
only if
\begin{equation}                            \label{idem}
v(\vartheta-L)\>+\>v(\vartheta-L)\>=\>v(\vartheta-L)\;.
\end{equation}
\end{theorem}
Note that $v(\vartheta-L)+v(\vartheta-L):=\{\alpha+\beta\mid\alpha,\beta
\in v(\vartheta-L)\}$ and that the sum of two initial segments of a
value group is again an initial segment. Equation (\ref{idem}) means
that $v(\vartheta-L)$ defines a cut in $vL$ which is idempotent under
addition of cuts (defined through addition of the left cut sets). If
$vL$ is archimedean, then there are only four possible idempotent cuts,
corresponding to $v(\vartheta-L)=\emptyset$ (which is impossible),
$v(\vartheta-L)=(vL)^{<0}$, $v(\vartheta-L)=(vL)^{\leq 0}$, and
$v(\vartheta-L)=vL$ (which means that $\vartheta$ lies in the completion
of $(L,v)$).

It is important to note that $v(\vartheta-K)\subseteq (vL)^{<0}$.
Indeed, if there were some $c\in K$ such that $v(\vartheta-c)\geq 0$,
then
\[
0\>\leq\>v((\vartheta-c)^p-(\vartheta-c))\>\leq\>
v(\vartheta^p-\vartheta-(c^p-c))\;.
\]
But a polynomial $X^p-X-a$ with $va\geq 0$ splits completely in the
absolute inertia field of $(L,v)$ and thus cannot induce a defect
extension. Therefore, if $vL$ is archimedean, then (\ref{idem}) holds if
and only if $v(\vartheta-L)=(vL)^{<0}$. This shows that the extension
$(L(\vartheta)|L,v)$ of Example~\ref{fse} is an independent
Artin-Schreier defect extension even if $L$ is not perfect. On the other
hand, the extension $(L(\vartheta_0)|L,v)$ of Example~\ref{transf},
where $\vartheta_0/b$ is a root of the polynomial $X^p-X-1/b^p t$, is a
dependent Artin-Schreier defect extension as it was obtained from the
purely inseparable defect extension $(L(\eta)|L,v)$. And indeed,
\begin{eqnarray*}
v\left(\frac{\vartheta_0}{b}\,-\,L\right)\,+\,v\left(\frac{\vartheta_0}{b}
\,-\,L\right) & = & \{\alpha\in vL\mid \alpha<vb\}\,+\,\{\alpha\in
vL\mid \alpha<vb\}\\
 & \ne & \{\alpha\in vL\mid \alpha<vb\}
\end{eqnarray*}
since $vb\ne 0$. Note that since $vL$ is $p$-divisible, we in fact have
that
\[\{\alpha\in vL\mid \alpha<vb\}+\{\alpha\in vL\mid \alpha<vb\}\>=\>
\{\alpha\in vL\mid \alpha<2vb\}\;.\]
This example also shows that the criterion of Theorem~\ref{dindcr} only
works for roots of Artin-Schreier polynomials. Indeed, $v(\vartheta_0-L)
= v(\eta-L)=(vL)^{<0}$, which does not contradict the theorem since the
minimal polynomial of $\vartheta_0$ is $X^p-dX-1/t$ with $vd\ne 0$.

\pars
Each of the perfect fields $(L,v)$ from Example~\ref{fse} provides an
example of a valued field without dependent Artin-Schreier defect
extensions, but admitting an independent Artin-Schreier defect
extension. Valued fields without independent Artin-Schreier defect
extensions but admitting dependent Artin-Schreier defect
extensions are harder to find; one example is given in \cite{[Ku9]}.

\parm
The classification of Artin-Schreier defect extensions and
Theorem~\ref{dindcr} are the main tool in the proof of the following
criterion (\cite{[Ku1]}; see \cite{[Ku2]} and \cite{[Ku9]}):

\begin{theorem}                             \label{dlamid}
A valued field $(L,v)$ of positive characteristic is henselian and
defectless if and only if it is algebraically maximal and inseparably
defectless.
\end{theorem}
This criterion is very useful when one tries to construct examples of
defectless fields with certain properties, as was done in
\cite[Section~4]{[Ku5]}.
How can one construct defectless fields? One possibility is to take any
valued field and pass to its maximal immediate extension. Every maximal
field is defectless. But it is in general an extension of very large
transcendence degree. If we want something smaller, then we could use
Theorem~\ref{ai}. But that talks only about function fields (or their
henselizations). If we want to construct a field with a certain value
group (as in \cite{[Ku5]}), we may have to pass to an infinite algebraic
extension. If we replace that by any of its maximal immediate algebraic
extensions, we obtain an algebraically maximal field $(M,v)$. But
Example~\ref{exdelon} shows that such a field may not be defectless.
If, however, we can make sure that $(M,v)$ is also inseparably
defectless, then Theorem~\ref{dlamid} tells us that $(M,v)$ is
defectless.

How do we know that a valued field $(L,v)$ is inseparably defectless?
In the case of finite \textbf{$p$-degree} $[K:K^p]$ (also called
\textbf{Ershov invariant of $K$}), Delon (\cite{[D]}) gave a
handy characterization of inseparably defectless valued fields:

\begin{theorem}                             \label{charinsdl}
Let $L$ be a field of characteristic $p>0$ and finite $p$-degree
$[L:L^p]$. Then for the valued field $(L,v)$, the property of being
inseparably defectless is equivalent to each of the following
properties:\sn
{\bf a)}\ \ $[L:L^p]=(vL:pvL)[Lv:Lv^p]$, i.e.,
$(L|L^p,v)$ is a defectless extension\sn
{\bf b)}\ \ $(L^{1/p}|L,v)$ is a defectless extension\sn
{\bf c)}\ \ every immediate extension of $(L,v)$ is separable\sn
{\bf d)}\ \ there is a separable maximal immediate extension of
$(L,v)$.
\end{theorem}

The very useful upward direction of the following lemma was also
stated by Delon (\cite{[D]}, Proposition 1.44):
\begin{lemma}                               \label{delon}
Let $(L'|L,v)$ be a finite extension of valued fields of characteristic
$p>0$. Then $(L,v)$ is inseparably defectless and of finite $p$-degree
if and only if $(L',v)$ is.
\end{lemma}
The condition of finite $p$-degree is necessary, as
Example~\ref{exampnagata} shows. In that example, $(k((t))|K,v)$ is a
purely inseparable defect extension of degree $p$. Hence $(K,v)$ is not
inseparably defectless. But $(k((t)),v)$ is, since it is a maximal and
hence defectless field.

\sn
{\bf Work in progress.}
\n
a) \ An analogue of the classification of Artin-Schreier defect extensions
and of Theorem~\ref{dlamid} has to be developed for the mixed
characteristic case (valued fields of characteristic 0 with residue
fields of characteristic $p$). Temkin and other authors have
already done part of the work.
\sn
b) \ The classification of Artin-Schreier defect extensions is also
reflected in their higher ramification groups. This will be worked out
in \cite{[Ku--Pi]}.

\parm
We have seen in Example~\ref{transf} that every purely inseparable
defect extension of degree $p$ of $(L,v)$ which does not lie in the
completion of $(L,v)$ can be transformed into a (dependent)
Artin-Schreier defect extension. This can be used to prove the following
result (cf.\ \cite{[Ku9]}):

\begin{proposition}
Assume that $(L,v)$ does not admit any dependent Artin-Schreier
extension. Then every immediate purely inseparable extension lies in the
completion of $(L,v)$.
\end{proposition}

\begin{corollary}
Every non-trivially valued Artin-Schreier closed field lies dense in its
perfect hull. In particular, the algebraic closure of a non-trivially
valued separable-algebraically closed field $(L,v)$ lies in the
completion of $(L,v)$.
\end{corollary}

\parm
Which of the Artin-Schreier defect extensions are the more harmful, the
dependent or the independent ones? There are some indications that the
dependent ones are more harmful. Temkin's work (especially
\cite[Theorem~3.2.3]{[T]}) seems to indicate that there is a
generalization of Theorem~\ref{hr} which already works over henselian
perfect instead of tame valued fields $(K,v)$. When $K$ is perfect, then
$(K,v)$ does not have dependent Artin-Schreier defect extensions. The
independent ones do not seem to matter here, at least when $(K,v)$ has
rank 1. This improvement is one of the keys to Theorem~\ref{lu3}. Let us
give an example which illustrates what is going on here.

\begin{example}
Assume that $(F|K,v)$ is an immediate function field of trans\-cendence
degree 1, rank 1 and characteristic $p>0$, and that we have chosen $x\in
F$ such that $(F^h|K(x)^h,v)$ is of degree $p$. We want to improve our
choice of $x$, that is, find $y\in F$ such that $F^h=K(y)^h$. The
procedure given in \cite{[Ku1],[Ku11]} uses the fact that because
$(F|K,v)$ is immediate, $x$ is a \textbf{pseudo limit} of a
\textbf{pseudo Cauchy sequence} $(a_\nu)_{\nu<\lambda}$ in $K$ which has
no pseudo limit in $K$ (\cite[Section~2]{[Ka]}). The hypothesis that
$(K,v)$ be tame (or separably tame and $F|K$ separable) guarantees that
$(a_\nu)_{\nu<\lambda}$ is of \textbf{transcendental type}, that is, if
$f$ is any polynomial in one variable over $K$, then the value
$vf(a_\nu)$ is fixed for all large enough $\nu<\lambda$
(\cite[p.~306]{[Ka]}). This is essential in our procedure. If we drop
the tameness hypothesis, then $(a_\nu)_{\nu<\lambda}$ may not be of
transcendental type and may in fact have some element in some immediate
algebraic extension of $(K,v)$ as a pseudo limit. Now suppose that this
element is the root $\vartheta$ of an Artin-Schreier polynomial over
$K$. The fact that both $x$ and $\vartheta$ are limits of
$(a_\nu)_{\nu<\lambda}$ implies that $v(x-\vartheta)>v(\vartheta-K)$. If
the Artin-Schreier defect extension $(K(\vartheta)|K,v)$ is independent,
then because of our rank 1 assumption, it follows that
$v(x-\vartheta)\geq 0$. If we assume in addition that $Kv$ is
algebraically closed, then there is some $c\in K$ such that
$v(x-\vartheta-c)> 0$. But then, by a special version of
Krasner's Lemma (cf.\ \cite[Lemma~2.21]{[Ku6]}), the polynomial
$X^p-X-(\vartheta^p-\vartheta)$ splits in $(K(x)^h,v)$, so that
$\vartheta\in K(x)^h$. This shows that $K$ is not relatively
algebraically closed in $K(x)^h$. Replacing $K$ by its relative
algebraic closure in $K(x)^h$, we will avoid this special case of pseudo
Cauchy sequences that are not of transcendental type.

If on the other hand, the extension $(K(\vartheta)|K,v)$ is dependent,
then it does not follow that $v(x-\vartheta)\geq 0$. But if
$v(x-\vartheta)<0$, Krasner's Lemma is of no use. However, by assuming
that $K$ is perfect we obtain that $(K,v)$ does not admit dependent
Artin-Schreier defect extensions. Assuming in addition that $K$ is
relatively algebraically closed in $F^h$, we obtain that
$(a_\nu)_{\nu<\lambda}$ does not admit any Artin-Schreier root
$\vartheta$ over $K$ as a limit. This fact alone does not imply
that under our additional assumptions, $(a_\nu)_{\nu<\lambda}$ must be
of transcendental type. But with some more technical effort, building on
results in \cite{[Ku9]}, this can be shown to be true.
\end{example}

\pars
Another indication may come from the paper \cite{[C--Pi]} by Steven Dale
Cutkosky and Olivier Piltant. They give an example of an extension
$R\subset S$ of algebraic regular local rings of dimension two over a
field $k$ of positive characteristic and a valuation on the rational
function field $\Quot R$, with $\Quot S|\Quot R$ being a tower of two
Artin-Schreier defect extensions, such that strong monomialization in
the sense of Theorem~4.8 of their paper does not hold for $R\subset S$
(\cite[Theorem~7.38]{[C--Pi]}).

Work in progress with Laura Ghezzi and Samar El-Hitti indicates that
both extensions are dependent Artin-Schreier defect extensions. In fact,
in Piltant's own words, he chose the valuation in the example such that
it is ``very close'' to [the behaviour of] a valuation in a purely
inseparable extension.

\sn
{\bf Open problem (CPE):} Is there a version of the example of Cutkosky
and Piltant involving independent Artin-Schreier defect extensions? Or
are such extensions indeed less harmful than the dependent ones? Can
strong monomialization always be proven when only independent
Artin-Schreier defect extensions are involved?

\end{section}

%
%
\begin{section}{Two languages}
Algebraic geometers and valuation theorists often speak different
languages. For example, while the defect was implicitly present already
in Abhyankar's work, it has been explicitly studied rather by the early
valuation theorists like Ostrowski, and by researchers interested in the
model theory of valued fields in positive characteristic or in constant
reduction, most of them using a field theoretic language and ``living in
the Kaplansky world of pseudo Cauchy sequences'' (cf.\
\cite[Section~2]{[Ka]}, \cite{[Ku2]}).

For instance, our joint investigation with Ghezzi and ElHitti of the
example given by Cutkosky and Piltant is facing the difficulty that
the valuation in the example is given by means of generating sequences,
whereas our criterion for dependence/independence is by nature closer to
the world of pseudo Cauchy sequences, which can also be used to describe
valuations.

\sn
{\bf Open problem (CGS):} Rewrite the criterion given in
Theorem~\ref{dindcr} in terms of generating sequences.

\pars
As to the problem of how to describe valuations, Michel Vaqui\'e has
done much work generalizing MacLane's approach using families of key
polynomials. In this approach, he also showed how to read off defects as
invariants of such families (see \cite{[V1]}). A closer look reveals
that a set like $v(\vartheta-L)$ can be directly determined from
Vaqui\'e's families of key polynomials.

\sn
{\bf Open problem (CV):} Is there an efficient algorithm to convert
generating sequences into Vaqui\'e's families of key polynomials? More
generally, find algorithms that convert between generating sequences,
key polynomials, pseudo Cauchy sequences and higher ramification groups.

\pars
A lot of work has been done by several authors on the description of
valuations on rational function fields, working with key polynomials or
pseudo Cauchy sequences. For references, see \cite{[Ku6]}.

\sn
{\bf Open problem (RFF):} Develop a thorough theory of valuations on
rational function fields, bringing the different approaches listed in
(CV) together, then generalize to algebraic function fields. Understand
the defect extensions that can appear over rational function fields.

\pars
Problems (CGV), (CV) and (RFF) can be understood as parts of a larger
program:
\sn
{\bf Open problem (DIC):} Develop a ``dictionary'' between algebraic
geometry and valuation theory. This would allow us to translate results
known about the defect into results in algebraic geometry, and open
questions in algebraic geometry into questions in valuation theory. It
would help us to investigate critical examples from several points of
view and to use them both in algebraic geometry and valuation theory.

\pars
Let us conclude with the following
\n
{\bf Open problem (DAG):} What exactly is the meaning of the defect in
algebraic geometry? How can we locate and interpret it? What is the
role of dependent and independent Artin-Schreier defect extensions,
e.g., in the work of Abhyankar, of of Cutkosky and Piltant, or of
Cossart and Piltant?

\end{section}

\end{document}